\documentclass[5p,times]{elsarticle}

\usepackage{lineno,hyperref}

\usepackage{graphicx}
\usepackage{amssymb}
\usepackage{tikz} 
\usepackage{amsmath}
\usepackage{algorithm, algorithmic}
\usepackage[font=normalsize]{caption} 
\usepackage{float}
\floatstyle{plaintop}
\restylefloat{table}

\journal{Parallel Computing}

\begin{document}

\begin{frontmatter}

\title{Scalable Algorithms for High-Order Approximations on Three-Dimensional Compact Stencils}

\author{Ronald Gonzales}
\author{Yury Gryazin\fnref{myfootnote}}
\author{Yun Teck Lee}
\address{Idaho State University, Department of Mathematics and Statistics, 921 S. 8th Ave., Stop 8085, Pocatello, ID 83209, USA}
\fntext[myfootnote]{Corresponding author, Email: gryazin@isu.edu}

\begin{abstract}
The recent development of multicore technologies on modern desktop computers makes parallelization of the proposed numerical approaches a priority in algorithmic research. The main performance improvement of personal computers in the upcoming years will be made based on the increasing number of cores on modern CPUs. This shifts the focus of algorithmic research from the development of sequential numerical methods to parallel methodology. 

This paper presents an efficient parallel direct algorithm with near-optimal complexity for the compact fourth and sixth-order approximation of the three-dimensional Helmholtz equations \citep{Turkel} with the problem coefficient depending on only one of the coordinate directions. The developed method is based on a combination of the separation of variables technique and a Fast Fourier Transform (FFT) type method. Similar direct solvers for the lower-order approximations of the two and three-dimensional Helmholtz equation were considered in several previous publications by the authors and other researchers (see e.g. \citep{gkl, yg2, eo, eo1}). 

The authors also consider a generalization of the presented algorithm to the solution of a wide class of linear systems obtained from approximation on the compact 27-point three-dimensional stencils on the rectangular grids with similar requirements on the stencil coefficients. The general restrictions on the coefficients in the considered class of compact schemes are developed and presented.
This class includes the second, fourth and sixth-order compact approximation schemes for the three-dimensional Helmholtz equation considered in this paper and our previous publications \citep{gkl, yg2, yg3}. As an example of the diversity of applications of the developed general method, the direct parallel implementation of a compact fourth-order approximation scheme for a convection-diffusion equation is considered. 

Another goal of this paper is to investigate the scalability of the proposed technique in the case of a large linear system using different parallel programming extensions.  The results of the implementation of this method in OpenMP, MPI and hybrid programming environments on the multicore computers and multiple node clusters are presented and discussed. The results demonstrate the high efficiency of the proposed direct solvers for many important applications on the structured grid with the corresponding 27-diagonal matrices of sizes up to $10^{11}$ by $10^{11}$.  \end{abstract}

\begin{keyword}
Compact finite-difference schemes \sep FFT \sep parallel algorithms \sep OpenMP \sep MPI \sep Hybrid
\end{keyword}

\end{frontmatter}


\section{Introduction}

In recent years, the problem of increasing the resolution of existing numerical solvers has become an urgent task in many areas of science and engineering. Most of the existing efficient solvers for structured matrices were developed for lower-order approximations of partial differential equations. The need for improved accuracy of the underlying algorithms leads to modified linear systems. As a result, the numerical solvers must be modified (see e.g. \citep{zs}). The time and memory constraints of practical applications make utilization of the existing sequential methods unacceptable in many situations. The development of the parallel algorithms for the implementation of the new high-resolution schemes becomes the necessary stage in the simulation process of many natural phenomena and engineering applications.

The focus of this paper is to introduce an efficient parallel direct solver for the recently developed sixth-order approximation compact schemes for the three-dimensional (3D) Helmholtz equation in both shared and distributed memory programming environments. For the two-dimensional (2D) constant-coefficient Helmholtz equation, this scheme was presented in \citep{nsd}, and the generalization to the 3D Helmholtz equation with constant coefficient was developed in \citep{sut}. The sixth-order schemes for 2D and 3D Helmholtz equations with nonconstant coefficients were proposed in \citep{Turkel} and the sixth-order approximation of Neumann and Sommerfeld-like boundary conditions were considered in \citep{yg2}. The second author considered the second-order approximation direct solver for 2D and 3D Helmholtz equations with the problem coefficient depending only on one coordinate direction as a preconditioner in general nonconstant coefficient problems in \citep{gkl, yg2}. 
The parallel direct algorithm for the second-order finite-difference and trilinear finite element approximations of the 3D constant-coefficient Helmholtz equation with Dirichlet and Neumann boundary conditions was considered in \citep{eo}. This algorithm was implemented using the MPI extension for Fortran90 and it demonstrated ``a large amount of parallelism'' \citep{eo}. The direct solver was used as a preconditioner for the problem with more general Sommerfeld-like boundary conditions. The average speedup in going from $\it{p}$ to $\it{2p}$ processors was between $1.77$ and $1.90$.  The results from \citep{eo, gkl, yg2} demonstrated the high efficiency of the FFT-based separation of variables approach and its high quality as a preconditioner in more general scattering problems. 

The first part of this paper extends this approach to the standard fourth-order compact scheme \citep{Lele92}, the sixth-order finite-difference approximation \citep{Turkel}, and proposes a parallel algorithm for the fast solution of the resulting systems. Then the authors generalize the developed parallel technique to the solution of a wider class of the linear systems obtained from approximations on 3D compact 27 point stencils with similar constraints on the stencil coefficients. This class includes the second, fourth and sixth-order compact approximation schemes for the three-dimensional Helmholtz equation considered in this paper and our previous publications \citep{gkl, yg2, yg3} as well as the trilinear finite element approximation considered in \citep{eo}. As an example of the diversity of applications of the developed general method, the direct parallel implementation of a compact fourth-order approximation scheme for a convection-diffusion equation is considered. The goal is to demonstrate the efficiency of the proposed direct parallel algorithms for the solution of a wide range of compact high-order approximation finite-difference and finite-elements schemes. The third and final goal of the paper is to demonstrate the advantages and limitations of the considered parallel approach in three popular parallel language extensions such as OpenMP, MPI and hybrid interfaces. 

The model problem considered in the paper is the numerical solution of
\begin{equation}
\nabla ^{2}u+k^{2}u=f , \ \  \textnormal{in}  \ \Omega, \label{problem}
\end{equation}
where $\Omega =\left\{ (x,y,z) \in \mathbb{R}^3 |L_x^l \leq x\leq L_x^r; L_y^l \leq y\leq L_y^r; L_z^l \leq z\leq L_z^r\right\},$ $L_x^l < L_x^r, L_y^l < L_y^r, L_z^l < L_z^r$ and $k$ is a complex valued coefficient depending only on $z$. The boundary conditions are
\begin{align}
&u=0, \ \textnormal {on}  \ \partial \Omega_1 = \left \{ (x,y,z) \in \Omega | \ x=L_x^l,L_x^r \ \textnormal {or} \ y=L_y^l,L_y^r  \right \}, \label{bc} \\
&\alpha_z \frac{\partial u}{\partial z} + \beta_z u=g(x,y), \ \textnormal {on} \ \partial \Omega_2 =\left \{ (x,y,z) \in \Omega | \ z=L_z^l,L_z^r \  \right \} , \nonumber   
\end{align}
where $\alpha_z, \beta_z\in \mathbb{C}$. These algorithms can be easily extended to the Neumann boundary conditions on the sides of the computational domain instead of the Dirichlet boundary conditions considered in this paper (see e.g. \citep{eo, eo1})

The resulting discretization leads to a system of linear equations with block 27-diagonal structure. In general, the matrix of this system is neither positive definite nor Hermitian. If an efficient direct solver can be applied then it is a natural choice. In the problem under consideration, the dependency of the problem coefficient on only one coordinate direction allows the use of the separation of variable technique based on FFT for an efficient direct solution of the discretized system. The well-known FFTW package \citep{FFTW_doc} for the standard implementation of the 2D discrete sine transform (DST) in the horizontal $xy-$ plane based on the FFT algorithm is utilized. This direct solver requires $O(N_xN_yN_z\log{N})$ operations, where $N = \max(N_x,N_y)$ and $N_x,N_y,$ and $N_z$ are the numbers of grid points in $x-$, $y-$, and $z-$ directions, respectively. 
 
The advantage of this approach is the natural parallelization of the proposed compact schemes on the DST and inverse DST steps. 
The transform implementation allows parallelization into the subdomains separated by the horizontal $xy-$ plane. It is possible to further parallelize this calculation by computing a set of 1D FFTs in the $x-$ direction then the $y-$ direction using a set of OpenMP threads or MPI processes. However, when this method was tested in OpenMP and MPI environments we found that the overhead drastically increased. Therefore, the first parallelization method was adapted. In the hybrid implementation of the solver on relatively large grids, the second method has some significant advantages. The approach alleviates a critical restriction on the number of MPI processes used by the direct solver in the distributed memory environment. So, in the hybrid algorithm, the second parallelization method was applied. On the step of the solution of the tridiagonal independent linear systems, parallelization is implemented by the division of the computational domain into a sequence of subdomains by the series of vertical $xz-$ and $yz-$planes. 

Numerical experiments with test problems demonstrate the high efficiency of the parallel implementation of the second, fourth and sixth-order compact finite-difference schemes. The performance of the developed methods is compared to the performance of the iterative block-parallel CARP-CG method \citep{Gordon} implemented on the multi-node cluster presented in \citep{Turkel}. The parallel algorithms presented here significantly outperformed the mentioned iterative method in both memory usage and ``wall time'' measure. 

The emphasis of the paper is on the parallel properties of the solver on the relatively large linear systems resulting from high-order 3D compact approximations. The largest considered grid in \citep{eo} was $60^3$ grid points, the authors in \citep{Turkel} restricted the consideration to grid sizes less than $402^3$. In contrast, this paper discusses the scalability of the presented algorithms on the 3D grids with up to $4096^3$ gridpoints. Particularly, the direct parallel double-precision solutions on the grid with $4096^3$ gridpoints was achieved in 28 sec on a 256-nodes cluster in the hybrid environment. The experimental results demonstrate the high scalability of the high-resolution parallel solver presented in this paper. The authors are not aware of similar scalability results on large grids for the 3D Helmholtz equation with a sufficiently large problem coefficient (sufficiently large frequency) that appeared in recent publications. 

This algorithm can be used in many important applications as a stand-alone parallel solver as well as a preconditioner in the more general solvers for the 3D Helmholtz equation with a nonconstant coefficient depending on three spatial variables and nonreflecting boundary conditions \citep{yg2, yg3}. The quality of the FFT-type preconditioners in scattering problems is addressed in \citep{eo1}. The application of the presented high-order approximation parallel solver as a preconditioner in the GMRES method for the general nonconstant-coefficient 3D Helmholtz equation with perfectly matched layer (PML) boundary conditions were presented by the authors at the SIAM Conference on Parallel Processing for Scientific Computing (PP20), February 12-15, 2020, Seattle, Washington, U.S. and will be submitted for publication in the near future.

The remainder of the paper is organized as follows. In Section \ref{com_disc_section} the description of the second, fourth and sixth-order approximation compact schemes are presented. Section \ref{FFT_section} focuses on the generalization of the developed parallel algorithm to the general compact stencil calculations with the required properties. Section \ref{parallel_section} presents the details of the OpenMP, MPI and Hybrid implementation of the developed algorithms. In Section \ref{results_section}, the effectiveness of the proposed parallel algorithms is demonstrated on a series of test problems.

\section{Compact Discretizations} \label{com_disc_section}

To introduce high-order approximation compact schemes for the solution of the three-dimensional Helmholtz equation (\ref{problem}) with the boundary conditions (\ref{bc}) we consider a grid  $ \Omega_h= \{ (x_i, y_j,z_k)\ |\  x_i=L_x^l +ih_x, y_j=L_y^l+jh_y, z_l=L_z^l+lh_z, i=1,...,N_x, j=1,...,N_y, l=1,...,N_z \}, $ where
 $ h_\alpha = \left (L_{\alpha}^r-L_{\alpha}^l \right )/(N_\alpha+1), \alpha=x,y,z $ are grid steps in $ \alpha-$ direction. The standard notation for the first and second-order central differences at $(i,j,l)-th$ grid point is given by
\begin{align}
\delta_{x}u_{i,j,l}=\frac{u_{i+1,j,l}-u_{i-1,j,l}}{2h_x}, \  \delta_x^2u_{i,j,l}=\frac{u_{i-1,j,l}-2u_{i,j,l}+u_{i+1,j,l}}{h_x^2}, \nonumber 
\end{align}
where $u_{i,j,l}=u(x_i,y_j,z_l)$.
The difference operators $\delta_{y}$, $\delta_{z}$, $\delta_{y}^2$ and $\delta_{z}^2$ used in the following sections are defined similarly. 
The second-order approximation scheme can be written as
\begin{align}
\left (\delta_x^2 + \delta_y^2 + \delta_z^2 \right)U_{i,j,l} +  k_{l}^2U_{i,j,l} = f_{i,j,l} \label{scheme2} 
\end{align}
where $k_{l} =k(z_l)$ and $U_{i,j,l}$ is the second-order finite-difference approximation to the solution $u_{i,j,l}$ of (\ref{problem}, \ref{bc}).

\subsection{Fourth-Order Pad\'e Approximation Compact Scheme} 

Now, we consider the standard fourth-order Pad\'e finite difference compact  approximation (see e.g. \cite{Lele92}) of (\ref{problem}). The fourth-order rational approximation of a second derivative can be presented in the form 
\begin{align}
\left . \frac{ \partial^2 u}{\partial \alpha^2} \right |_{i,j,l} =\left (1+\frac{h_\alpha^2}{12} \delta_\alpha^2 \right)^{-1} \delta_\alpha^2u_{i,j,l}  + O\left(h_\alpha^4\right), \ \alpha=x,y,z. \nonumber
\end{align}
By substituting this approximation into (\ref{problem}), we obtain
\begin{align*}
&\left (1+\frac{h_x^2}{12} \delta_x^2 \right)^{-1} \delta_x^2u_{i,j,l} +\left (1+\frac{h_y^2}{12} \delta_y^2 \right)^{-1} \delta_y^2u_{i,j,l} +  \\ 
&\left (1+\frac{h_z^2}{12} \delta_z^2 \right)^{-1} \delta_z^2u_{i,j,l} + k_{l}^2u_{i,j,l} = f_{i,j,l} + O\left(\max\left(h_x^4,h_y^4,h_z^4\right)\right),
\end{align*}
which gives
\begin{align}
&\left (\delta_x^2 + \delta_y^2 + \delta_z^2 \right)U_{i,j,l} + \frac{(h_x^2 +h_y^2)}{12}\delta_x^2\delta_y^2 U_{i,j,l}+ \frac{(h_x^2 +h_z^2)}{12}\delta_x^2\delta_z^2U_{i,j,l} \nonumber\\
&+ \frac{(h_y^2 +h_z^2)}{12}\delta_y^2\delta_z^2 U_{i,j,l} + \left (1+\frac{h_x^2}{12} \delta_x^2 + \frac{h_y^2}{12} \delta_y^2 + \frac{h_z^2}{12} \delta_z^2 \right) ( k_{l}^2U_{i,j,l}) \nonumber \\
&=\left (1+\frac{h_x^2}{12} \delta_x^2 + \frac{h_y^2}{12} \delta_y^2 + \frac{h_z^2}{12} \delta_z^2 \right)f_{i,j,l}=f^{(IV)}_{i,j,l} \label{scheme4}
\end{align}
where $U_{i,j,l}$ is the fourth-order compact finite-difference approximation to $u_{i,j,l}$.

\subsubsection{Sixth-Order Approximation Compact Scheme}

In this section we present a 3D compact sixth-order approximation finite-difference scheme. The scheme requires a uniform grid step so we assume $ h=h_x=h_y=h_z$. Using the appropriate derivatives of (\ref{problem}) we can write the sixth-order compact approximation of the equation in the form
\begin{align}
&\left( \delta_x^2 + \delta_y^2 + \delta_z^2\right)\left(1+ \frac{k_{l}^2h^2}{30}\right)U_{i,j,l}+ \frac{h^4}{30}\delta_x^2\delta_y^2\delta_z^2U_{i,j,l}+k_{l}^2U_{i,j,l} + \nonumber \\
&\frac{h^2}{6}( \delta_x^2\delta_y^2+ \delta_x^2\delta_z^2+ \delta_y^2\delta_z^2)\left (1+ \frac{k_{l}^2h^2}{15}\right )U_{i,j,l}+ \frac{h^2}{20}\Delta_h(k_{l}^2U_{i,j,l})= \nonumber \\
&f_{i,j,l}+ \frac{h^2}{12}\nabla^2 f_{i,j,l}+\frac{h^4}{360}\nabla^4 f_{i,j,l} +  \nonumber \\
&\frac{h^4}{90}\left ( \frac{\partial^4f}{\partial x^2\partial y^2}+\frac{\partial^4f}{\partial x^2\partial z^2}+\frac{\partial^4f}{\partial y^2\partial z^2}\right )_{i,j,l} = f^{(VI)}_{i,j,l} , \label{scheme6}
\end{align}
\noindent where
\begin{align}
&\Delta_h\left(k_{l}^2U_{i,j,l}\right) = k_{l}^2U_{i,j,l} + \left ( \frac{\partial^2}{\partial z^2} (k^2) -k^4 \right)_{l}  + \nonumber\\
&2[(k^2)_z]_{l} \left [ \delta_z \left ( 1 + \frac{h^2}{6} ( \delta_x^2 + \delta_y^2  + k^2_{l}) \right ) U_{i,j,l} -   \frac{h^2}{6}(f_z)_{i,j,l} \right ] . \nonumber
\end{align}

This compact scheme was developed in  \citep{Turkel} for the approximation of the 3D Helmholtz equation with nonconstant coefficient and Dirichlet boundary conditions. The compact sixth-order approximation of the Sommerfeld-like boundary conditions was developed in our previous publication \citep{yg2}. The partial derivatives of $f(x,y,z)$ in (\ref{scheme6}) are approximated by using the implicit compact approximation technique described in \citep{Lele92}. It requires boundary conditions on the corresponding derivatives. In many important applications, we can impose zero boundary conditions on these derivatives.

\section{Scalable FFT Compact Direct Solver} \label{FFT_section}

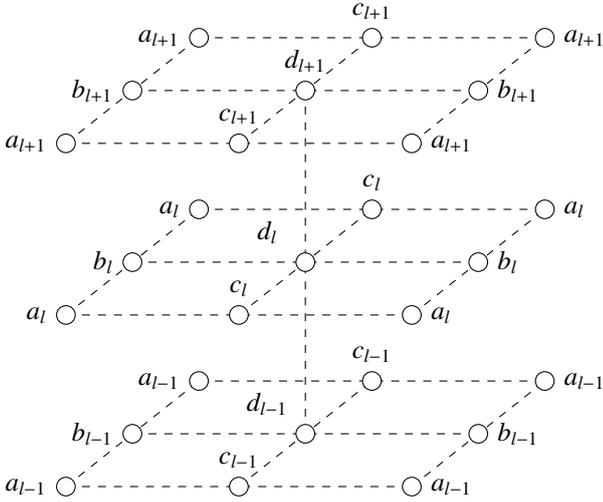
\begin{figure}[H]
\centering
\caption{27-Point Stencil}

\begin{tikzpicture}[xscale=1.75,yscale=1.75]

\def\length{1.3};\def\offsetx{.5};\def\offsety{.4};\def\vertoffset{-1.3};\def\dotsize{7pt};
\coordinate (top_1) at (0,0);
\coordinate (top_2) at (\length,0);
\coordinate (top_3) at (2*\length,0);
\coordinate (top_4) at (\offsetx,\offsety);
\coordinate (top_5) at (\length+\offsetx,\offsety);
\coordinate (top_6) at (2*\length+\offsetx,\offsety);
\coordinate (top_7) at (2*\offsetx,2*\offsety);
\coordinate (top_8) at (\length+2*\offsetx,2*\offsety);
\coordinate (top_9) at (2*\length+2*\offsetx,2*\offsety);

\coordinate (mid_1) at (0,\vertoffset);
\coordinate (mid_2) at (\length,\vertoffset);
\coordinate (mid_3) at (2*\length,\vertoffset);
\coordinate (mid_4) at (\offsetx,\offsety+\vertoffset);
\coordinate (mid_5) at (\length+\offsetx,\offsety+\vertoffset);
\coordinate (mid_6) at (2*\length+\offsetx,\offsety+\vertoffset);
\coordinate (mid_7) at (2*\offsetx,2*\offsety+\vertoffset);
\coordinate (mid_8) at (\length+2*\offsetx,2*\offsety+\vertoffset);
\coordinate (mid_9) at (2*\length+2*\offsetx,2*\offsety+\vertoffset);

\coordinate (bot_1) at (0,2*\vertoffset);
\coordinate (bot_2) at (\length,2*\vertoffset);
\coordinate (bot_3) at (2*\length,2*\vertoffset);
\coordinate (bot_4) at (\offsetx,\offsety+2*\vertoffset);
\coordinate (bot_5) at (\length+\offsetx,\offsety+2*\vertoffset);
\coordinate (bot_6) at (2*\length+\offsetx,\offsety+2*\vertoffset);
\coordinate (bot_7) at (2*\offsetx,2*\offsety+2*\vertoffset);
\coordinate (bot_8) at (\length+2*\offsetx,2*\offsety+2*\vertoffset);
\coordinate (bot_9) at (2*\length+2*\offsetx,2*\offsety+2*\vertoffset);

\draw[dashed] (top_1)--(top_7);
\draw[dashed] (top_2)--(top_8);
\draw[dashed] (top_3)--(top_9);
\draw[dashed] (top_1)--(top_3);
\draw[dashed] (top_4)--(top_6);
\draw[dashed] (top_7)--(top_9);

\draw[dashed] (mid_1)--(mid_7);
\draw[dashed] (mid_2)--(mid_8);
\draw[dashed] (mid_3)--(mid_9);
\draw[dashed] (mid_1)--(mid_3);
\draw[dashed] (mid_4)--(mid_6);
\draw[dashed] (mid_7)--(mid_9);

\draw[dashed] (bot_1)--(bot_7);
\draw[dashed] (bot_2)--(bot_8);
\draw[dashed] (bot_3)--(bot_9);
\draw[dashed] (bot_1)--(bot_3);
\draw[dashed] (bot_4)--(bot_6);
\draw[dashed] (bot_7)--(bot_9);

\draw[dashed] (top_5)--(bot_5);

   \node[circle,draw=black, fill=white, inner sep=0pt,minimum size=\dotsize,label=left:{$a_{l+1}$}] () at (top_1) {};
\node[circle,draw=black, fill=white, inner sep=0pt,minimum size=\dotsize,label=above:{$c_{l+1}$}] () at (top_2) {};  
 \node[circle,draw=black, fill=white, inner sep=0pt,minimum size=\dotsize,label=right:{$a_{l+1}$}] () at (top_3) {};       
  \node[circle,draw=black, fill=white, inner sep=0pt,minimum size=\dotsize,label=left:{$b_{l+1}$}] () at (top_4) {};  
    \node[circle,draw=black, fill=white, inner sep=0pt,minimum size=\dotsize,label=above:{$d_{l+1}$}] () at (top_5) {};  
        \node[circle,draw=black, fill=white, inner sep=0pt,minimum size=\dotsize,label=right:{$b_{l+1}$}] () at (top_6) {};  
    \node[circle,draw=black, fill=white, inner sep=0pt,minimum size=\dotsize,label=left:{$a_{l+1}$}] () at (top_7) {};  
        \node[circle,draw=black, fill=white, inner sep=0pt,minimum size=\dotsize,label=above:{$c_{l+1}$}] () at (top_8) {};     
 \node[circle,draw=black, fill=white, inner sep=0pt,minimum size=\dotsize,label=right:{$a_{l+1}$}] () at (top_9) {};

\node[circle,draw=black, fill=white, inner sep=0pt,minimum size=\dotsize,label=left:{$a_{l}$}] () at (mid_1) {};
\node[circle,draw=black, fill=white, inner sep=0pt,minimum size=\dotsize,label=above:{$c_{l}$}] () at (mid_2) {};  
 \node[circle,draw=black, fill=white, inner sep=0pt,minimum size=\dotsize,label=right:{$a_{l}$}] () at (mid_3) {};       
  \node[circle,draw=black, fill=white, inner sep=0pt,minimum size=\dotsize,label=left:{$b_{l}$}] () at (mid_4) {};  
    \node[circle,draw=black, fill=white, inner sep=0pt,minimum size=\dotsize,label={[shift={(-.5,0)}]$d_{l}$}] () at (mid_5) {};  
        \node[circle,draw=black, fill=white, inner sep=0pt,minimum size=\dotsize,label=right:{$b_{l}$}] () at (mid_6) {};  
    \node[circle,draw=black, fill=white, inner sep=0pt,minimum size=\dotsize,label=left:{$a_{l}$}] () at (mid_7) {};  
        \node[circle,draw=black, fill=white, inner sep=0pt,minimum size=\dotsize,label=above:{$c_{l}$}] () at (mid_8) {};     
 \node[circle,draw=black, fill=white, inner sep=0pt,minimum size=\dotsize,label=right:{$a_{l}$}] () at (mid_9) {};

\node[circle,draw=black, fill=white, inner sep=0pt,minimum size=\dotsize,label=left:{$a_{l-1}$}] () at (bot_1) {};
\node[circle,draw=black, fill=white, inner sep=0pt,minimum size=\dotsize,label=above:{$c_{l-1}$}] () at (bot_2) {};  
 \node[circle,draw=black, fill=white, inner sep=0pt,minimum size=\dotsize,label=right:{$a_{l-1}$}] () at (bot_3) {};       
  \node[circle,draw=black, fill=white, inner sep=0pt,minimum size=\dotsize,label=left:{$b_{l-1}$}] () at (bot_4) {};  
    \node[circle,draw=black, fill=white, inner sep=0pt,minimum size=\dotsize,label={[shift={(-.5,0)}]$d_{l-1}$}] () at (bot_5) {};  
        \node[circle,draw=black, fill=white, inner sep=0pt,minimum size=\dotsize,label=right:{$b_{l-1}$}] () at (bot_6) {};  
    \node[circle,draw=black, fill=white, inner sep=0pt,minimum size=\dotsize,label=left:{$a_{l-1}$}] () at (bot_7) {};  
        \node[circle,draw=black, fill=white, inner sep=0pt,minimum size=\dotsize,label=above:{$c_{l-1}$}] () at (bot_8) {};     
 \node[circle,draw=black, fill=white, inner sep=0pt,minimum size=\dotsize,label=right:{$a_{l-1}$}] () at (bot_9) {};

\end{tikzpicture}

\label{27_point_stencil}
\end{figure}

In this section, we present a generalized parallel direct solver for the second, fourth and sixth-order compact schemes presented in Section 2. This algorithm was developed for the solution of a 27-diagonal linear system satisfying a set of required conditions. Figure \ref{27_point_stencil} demonstrates the 3D compact stencil form corresponding to this system.

\subsection{Stencil Form of the Schemes}
All three compact schemes under consideration can be presented in the displayed stencil form. Also, any compact scheme with these stencil coefficients can be expressed at every grid point $(i,j,l)$ as 
\begin{align}
\Sigma^{\nu=l+1}_{\nu=l-1}  & \left( a_{\nu}\left[ U_{i-1,j-1,\nu} + U_{i-1,j+1,\nu} + U_{i+1,j-1,\nu} + U_{i+1,j+1,\nu}\right] + \right . \nonumber \\
 & \ \ b_{\nu}\left[ U_{i-1,j,\nu} + U_{i+1,j-1,\nu} \right] + c_{\nu}\left[ U_{i,j-1,\nu} + U_{i,j+1,\nu} \right]  \label{stencil27} \\
 & \ \  \left. d_{\nu} U_{i,j,\nu} \right) = F_{i,j,l}.\nonumber
\end{align}
This equation corresponds to the $(i+N_x \cdot j + N_x \cdot N_y \cdot l)-th$ row in the resulting linear system $AU=F.$
The following subsections specify the stencil coefficients in all three compact schemes under consideration.

\subsection{Second-Order}

Let $R_{zx} = h_z^2/h_x^2$ and $R_{zy} = h_z^2/h_y^2$. In the second-order scheme, there are only five nonzero parameters
\begin{align*}
b_{l}=R_{zx} , c_{l} =R_{zy}, d_{l-1} = d_{l+1}=1, d_l = -2\left(R_{zx}+ R_{zy} +1\right) + h_z^2k^2_l.
\end{align*}
\noindent{The} right hand side is given by $F_{i,j,l}=h_z^2f_{i,j,l}.$

\subsection{Fourth-Order}

\noindent{The} nonzero coefficients in the fourth-order scheme are 
\setlength{\mathindent}{0pt}
\begin{align*}
&b_{l-1} = b_{l+1} = (1 + R_{zx})/12, \ \ c_{l-1} = c_{l+1} = (1 + R_{zy})/12,\\
&d_{l-1} = 2/3 - ( R_{zx} + R_{zy})/6 + h_z^2k^2_{l-1}/12,\\
&d_{l+1} = 2/3 - ( R_{zx} + R_{zy})/6 + h_z^2k^2_{l+1}/12, \\ 
&a_{l} = (R_{zx}+R_{zy})/12, \ \ b_{l} = (4R_{zx} - R_{zy} - 1 + h_z^2k^2_l/2)/6,\\
&c_{l}=(4R_{zy} - R_{zx} - 1 + h_z^2k^2_l/2)/6, \\ 
&d_{l} =- 4(1 + R_{zx} + R_{zy})/3 + h_z^2k^2_l/2.
\end{align*}
Then the fourth-order right hand side is given by $F_{i,j,l}=h_z^2f^{(IV)}_{i,j,l}.$

\subsection{Sixth-Order}

\noindent{In} the case of the sixth-order scheme (\ref{scheme6}), the grid step size $h$ is uniform in all three directions. The scheme parameters are calculated as
\begin{align*}
&a_{l-1} = a_{l+1} = 1/30, \\
&b_{l-1} = c_{l-1} =1/10+ h^2k^2_{l-1}/90- h^3k^2_{z,l}/120,\\
&b_{l+1} = c_{l+1} = 1/10+ h^2k^2_{l+1}/90+h^3k^2_{z,l}/120,\\
&d_{l-1}=7/15-h^2k^2_{l-1}/90-  \left(h^3k^2_{z,l}/20\right)\left(1/3 + h^2k^2_{l-1}/6\right),\\
&d_{l+1}=7/15-h^2k^2_{l+1}/90+ \left(h^3k^2_{z,l}/20\right)\left(1/3 + h^2k^2_{l+1}/6\right),\\
&a_{l}  = 1/10+h^2k^2_l/90, \ \ b_{l}=c_{l} = 7/15-h^2k^2_l/90,\\
&d_{l} = -64/15 + 14h^2k^2_l/15 - h^4(k^2_l)^2/20 + h^4k^2_{zz,l}/20.
\end{align*}
The right hand side is given by $F_{i,j,l}=h_z^2f^{(VI)}_{i,j,l}.$

The stencil pattern (\ref{stencil27}) is not restricted to the algorithms presented in this paper. The trilinear finite element discretization of the weak formulation of the 3D Helmholtz equation with constant coefficient on rectangular grids (see e.g. \citep{eo}) can be presented in the same compact stencil form with the corresponding coefficient presented in (p.169, \citep{eo}). 

\subsection{Fourth-Order Approximation Scheme for 3D Convection-Diffusion Equation}

In this subsection, the versatility of the proposed parallel solver is illustrated on the 3D convection-diffusion equation with dominant convection in the $z-$ direction.
The steady-state 3D convection-diffusion equation can be written as
\begin{equation}
\nabla^2 u + \alpha \frac{\partial u}{\partial x} +\beta\frac{\partial u}{\partial y} +\gamma \frac{\partial u}{\partial z} = f(x,y,z), \label{cd_eq} \\
\end{equation}
where $\alpha$, $\beta$ and $\gamma$ are variable or constant convection coefficients in the $x-$, $y-$ and $z-$ directions respectively and $f$ is a forcing function. We assume that the horizontal gradient is significantly smaller than the first derivative of $u$ in the $z-$ direction. This is a common situation in the modeling of atmospheric heat convection.
Then the equation (\ref{cd_eq}) becomes  
\begin{equation}
\nabla^2 u + \gamma \frac{\partial u}{\partial z} = f(x,y,z). \label{cd_eq1} \\
\end{equation}

Now we extend the work done in \cite{conv_diff} on the 2D convection diffusion equation to the 3D case. Using the relevant derivatives of the original equation (\ref{cd_eq}), the compact fourth-order approximation scheme can be presented in the form (\ref{stencil27}) with the following nonzero stencil coefficients
\begin{align*}
&b_{l\pm1} = (1 + R_{zx})(2 \pm \gamma h_z)/24, \ \ c_{l\pm1} = (1 + R_{zy})(2 \pm \gamma h_z)/24,\\
&d_{l\pm1} = 2/3 - ( R_{zx} + R_{zy})/6 \pm \frac{h_z\gamma}{12}(4 - R_{zx} - R_{zy} \pm h_z\gamma)),\\
&a_{l} = (R_{zx}+R_{zy})/12, \ \ b_{l} = (4R_{zx} - R_{zy} - 1)/6,\\
&c_{l}=(4R_{zy} - R_{zx} - 1)/6, d_{l} =- 4(1 + R_{zx} + R_{zy})/3 - h_z^2\gamma^2/6.
\end{align*}
Then the right hand side of the resulting linear system can be presented as
\begin{equation}
F_{i,j,l}=h_z^2 \left [ f + \frac{h_x^2}{12} \frac{\partial^2 f}{\partial x^2} + \frac{h_y^2}{12} \frac{\partial^2 f}{\partial y^2}
+\frac{h_z^2}{12}\left(\gamma \frac{\partial f}{\partial z}+ \frac{\partial^2 f}{\partial z^2}\right) \right ]_{i,j,k} \label{cd_rhs}. \\
\end{equation}
For simplicity, the right hand side of (\ref{cd_eq1}) is assumed to be twice continuously differentiable. 

\subsection{FFT Solver for Compact Stencil}

The following derivation presents an efficient way to parallelize the proposed direct solver. The numerical scheme (\ref{stencil27}) can be presented in block three diagonal form written as
\begin{align*}
&C_1U_1+C_{p,1}U_{2}=F_1,\\ 
&C_{m,l}U_{l-1}+C_lU_l+C_{p,l}U_{l+1}=F_l , l=2,N_z-1,\\
&C_{m,N_z}U_{N_z-1}+C_{N_z}U_{N_z}=F_{N_z}.
\end{align*}
Here, the vectors $U_{l}$ and $F_{l}$ are the sections of the unknown vector $U$ and the right hand side $F$ with $l=1,...,N_z$.
The nine diagonal matrices $C_{m,l}, C_l,$ and $C_{p,l}$ are defined by the coefficients in (\ref{stencil27}). These matrices can be simultaneously diagonalized by using the $N_x \cdot N_y \times N_x \cdot N_y$ orthogonal matrix of eigenvectors $V$ defined by
\[V_{i,j}^{n,m}=\frac{2}{\sqrt{(N_x+1)(N_y+1)}}\sin\left(\frac{\pi n i}{N_x+1}\right)\sin\left(\frac{\pi m j}{N_y+1}\right),\] 
where  $1\le i,n\le N_x$ and $1\le j,m\le N_y$. The corresponding eigenvalues $\lambda_{i,j,\nu}$ for the matrices $C_{m,l}, C_l,$ and $C_{p,l}$ are given by
\begin{align*}
&\lambda_{i,j,\nu} = 4a_{\nu}\cos\left(\frac{(i+1)\pi}{N_x+1}\right)\cos\left(\frac{(j+1)\pi}{N_y+1}\right) + \\
&2b_{\nu}\cos\left(\frac{(i+1)\pi}{N_x+1}\right) + 2c_{\nu}\cos\left(\frac{(j+1)\pi}{N_y+1}\right)+d_{\nu}, \\
&0\le i < N_x, \ \ 0\le j< N_y, \ \ \nu = l-1,l,l+1.
\end{align*}
Since the matrices $\Lambda_m= V^TC_mV$, $\Lambda= V^TCV$ and $\Lambda_p= V^TC_pV$ are the diagonal matrices of eigenvalues, the original system can be presented as a set of $N_x \cdot N_y$ independent linear systems of size $N_z$ by $N_z$ in the following manner
\begin{align*}
&C_mU_{l-1}+CU_l+C_pU_{l+1}=F_l, \\
&V^TC_mVV^TU_{l-1}+V^TCVV^TU_l+V^TC_pVV^TU_{l+1}=V^TF_l, \\
&\Lambda_mW_{l-1}+\Lambda W_l +\Lambda_pW_{l+1} = \overline{F_l},
\end{align*}
where $W_l=V^TU_l$ and $\overline{F_l} =V^T F_l$. 

Each independent system in the set is tridiagonal and can be efficiently solved using LU decomposition with $O(N_z)$ computational complexity. The solution of each system in the set is independent of each other. Therefore it can be efficiently parallelized on multicore CPUs and clusters.

Prior to solving these independent systems, the transformed right-hand side vectors $\overline{F}_l=V^TF_l, l=1,..., N_z$ must be found. The matrix-vector multiplication in this calculation can be seen as a 2D DST of the right-hand side vector $F_l$.  This transform can be found by using the FFT algorithm with computational complexity $O(N_x \cdot N_y \textrm{ln}(N_x \cdot N_y))$. In our solver, we used the standard implementation of FFT from the open-source C library developed at Massachusetts Institute of Technology, namely FFTW \cite{FFTW_doc}.

\section{Parallelization}\label{parallel_section}

In this section, the details of the OpenMP, MPI and Hybrid implementations of the developed direct solvers are considered. We discuss the limitations and advantages of each implementation depending on the particular computer architecture. The goal is to demonstrate the scalability of the developed methods on modern multicore desktops and multi-node clusters.

\subsection{OpenMP}

First, we consider the parallelization of the direct solver using OpenMP, an application programming interface (API). OpenMP makes use of shared memory architecture and thus allows every thread to access all allocated memory in the program. While this is a convenient parallelization tool, it is restrictive. Programs using strictly OpenMP can only run on a single computer with shared memory. On a large multi-node cluster this typically restricts the parallel execution to a single node with 16 to 32 processors. However, a more significant limitation is the amount of random access memory (RAM) available on a single machine or a node. In many situations, the computations require vast amounts of RAM that are simply unavailable for OpenMP applications.

Despite its limitations, the shared address model allows a relatively simple implementation and excellent speed up in the execution of structured blocks. Algorithm \ref{3D_omp_alg} shows how OpenMP was used to implement the developed parallel direct solver. We can see that the algorithm is naturally divided into three easily parallelizable sections: forward DST step, the solution of the set of independent tridiagonal systems and the inverse DST step. In all three steps, the OpenMP threads use the different parts of the shared arrays. A rearrangement of the working arrays between the three consecutive stages of the algorithm may take significant processing time and must be implemented with careful consideration of array distribution between different types of CPU memory. 

\begin{algorithm}[H]
\caption{OpenMP 3D Helmholtz Direct Solver}
\begin{algorithmic} [1]\label{3D_omp_alg}
\STATE \#pragma omp parallel for
\FOR{$l=1,\dots,N_z$}
\STATE 2D forward DST in $x-,y-$ direction
\ENDFOR
\STATE \#pragma omp parallel for collapse(2)
\FOR{$j=1,\dots,N_y;i=1,\dots,N_x$}
\STATE Solve the tridiagonal system using LU decomposition
\ENDFOR
\STATE \#pragma omp parallel for
\FOR{$l=1,\dots,N_z$}
\STATE 2D inverse DST in $x-,y-$ direction
\ENDFOR
\end{algorithmic} 
\end{algorithm}

In the OpenMP implementation of the developed direct algorithm, a minor complication emerged. An FFTW plan is a necessary function that sets up the calculation of the FFT \cite{FFTW_doc} used in DST forward and inverse steps. These plans are not ``thread-safe'' therefore must be created within a critical region in the parallel section of the code.

Overall, the OpenMP implementation on a single desktop computer with a multicore CPU or a single multicore node of a cluster demonstrates excellent, near-linear scalability. It is perfect for a medium-sized grid. We were able to run our test problems with computational grids up to  $512^3$ on a single machine with 16G RAM. 

\subsection{MPI} \label{sec: MPI}

OpenMP provides a very convenient and efficient standard for parallel programming in the shared-memory environment, but for a large enough computational grid, the memory required to allocate the necessary arrays can overrun the RAM available on a single node. The natural solution to this problem is to distribute the working arrays and computational tasks between the nodes of a cluster.  In this case, there is no way for one processor to directly access the address space of another. This requires explicit message passing (MP), i.e.\,communication between processors. Several APIs were developed for this, but the standard today is the MPI (Message Passing Interface). 

The developed parallel algorithm is well suited for this type of parallelization since the different parallel processes are using different parts of the computational arrays during the program execution. In the developed MPI implementation of the algorithm, the sequential program was modified to run on several nodes allocating only the minimum required memory on each. This was accomplished by dividing the computational domain as evenly as possible along the vertical direction on the DST steps and in the horizontal directions on the independent tridiagonal solver steps. In turn, this enables the use of much larger computational grids as the program is no longer limited by the memory of a single node. The usual limitation in this implementation is the communication time. Each process runs the entire program on the assigned (distributed) section of the available memory independently and communicates with other processes only when explicitly specified. 
\begin{algorithm}[H]
\caption{MPI 3D Helmholtz Direct Solver}
\begin{algorithmic} [1] \label{MPI}
\STATE Find $start_y,start_z,end_y,$ and $end_z$ using the rank
\FOR{$l=start_z,\dots,end_z$}
\STATE 2D forward DST in $x-,y-$ direction
\ENDFOR
\STATE Scatter the data via MPI to the appropriate process
\FOR{$ j=start_y,\dots,end_y; i=1,\dots,N_x$}
\STATE Solve the tridiagonal system using LU decomposition
\ENDFOR
\STATE Scatter the data via MPI to the appropriate process
\FOR{$l=start_z,\dots,end_z$}
\STATE 2D inverse DST in $x-,y-$ direction
\ENDFOR
\end{algorithmic}
\end{algorithm}

 The communication portion of the program runtime (wall time) grows with the number of nodes. This becomes a major obstacle to the linear scalability of the designed algorithm. Algorithm \ref{MPI} shows how MPI was used to parallelize the direct solver. 

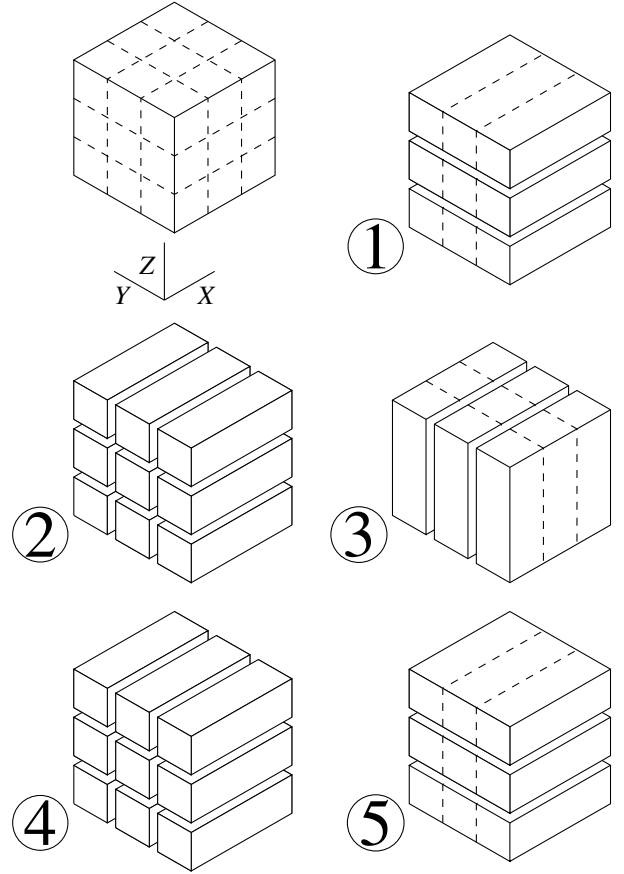
\begin{figure}[H]
\centering
\caption{Data Transfer between MPI Processes}

\begin{center}
\begin{tikzpicture}[scale=\textwidth/36cm,y={(0.866cm,-0.5cm)},x={(0.866cm,0.5cm)},z={(0cm,1cm)}]

\def\bigfont{4}
\def\h{1}
\def\w{3*\h}
\def\buff{.25}
\def\os{\h+\buff}
\def\ol{-.5}


\def\ox{-5}
\def\oy{-5}
\def\oz{1.35}

\def\a{\h*1}
\def\xoffset{-.15}
\def\yoffset{-.15}
\def\zoffset{-\a*1.75}
\coordinate (c1) at (\ox+\xoffset,\oy+\w+\yoffset,\oz+\zoffset);
\coordinate (c2) at (\ox+\h/2+\a+\xoffset,\oy+\w+\yoffset,\oz+\zoffset);
\coordinate (c3) at (\ox+\xoffset,\oy+\w-\a-\h/2+\yoffset,\oz+\zoffset);
\coordinate (c4) at (\ox+\xoffset,\oy+\w+\yoffset,\oz+\zoffset+\a+\h/2);
\draw[-] (c1)--(c2);
\draw[-] (c1)--(c3);
\draw[-] (c1)--(c4);
\node at (\ox+\xoffset,\oy+\yoffset+\w-\h/2,\oz+\zoffset+\h/1.5) {$Z$};
\node at (\ox+\xoffset-\h/2,\oy+\yoffset+\w-\h/4-\a/2,\oz+\zoffset) {$Y$};
\node at (\ox+\xoffset+\h/4+\a/2,\oy+\yoffset+\w+\h/2,\oz+\zoffset) {$X$};

\coordinate (b1) at (\ox,\oy,\oz+\w);
\coordinate (b2) at (\ox,\oy,\oz);
\coordinate (b3) at (\ox,\oy+\w,\oz);
\coordinate (b4) at (\ox,\oy+\w,\oz+\w);
\coordinate (b5) at (\ox+\w,\oy,\oz+\w);
\coordinate (b6) at (\ox+\w,\oy+\w,\oz+\w);
\coordinate (b7) at (\ox+\w,\oy+\w,\oz);
\draw[-] (b1)--(b2)--(b3)--(b4)--(b1)--(b5)--(b6)--(b7)--(b3);
\draw[-] (b4)--(b6);
\coordinate (c1) at (\ox,\oy+\h,\oz);
\coordinate (c2) at (\ox,\oy+\h,\oz+\w);
\coordinate (c3) at (\ox+\w,\oy+\h,\oz+\w);
\draw[dashed] (c1)--(c2)--(c3);
\coordinate (c1) at (\ox,\oy+\h*2,\oz);
\coordinate (c2) at (\ox,\oy+\h*2,\oz+\w);
\coordinate (c3) at (\ox+\w,\oy+\h*2,\oz+\w);
\draw[dashed] (c1)--(c2)--(c3);
\coordinate (c1) at (\ox,\oy,\oz+\h);
\coordinate (c2) at (\ox,\oy+\w,\oz+\h);
\coordinate (c3) at (\ox+\w,\oy+\w,\oz+\h);
\draw[dashed] (c1)--(c2)--(c3);
\coordinate (c1) at (\ox,\oy,\oz+\h*2);
\coordinate (c2) at (\ox,\oy+\w,\oz+\h*2);
\coordinate (c3) at (\ox+\w,\oy+\w,\oz+\h*2);
\draw[dashed] (c1)--(c2)--(c3);
\coordinate (c1) at (\ox+\h,\oy,\oz+\w);
\coordinate (c2) at (\ox+\h,\oy+\w,\oz+\w);
\coordinate (c3) at (\ox+\h,\oy+\w,\oz);
\draw[dashed] (c1)--(c2)--(c3);
\coordinate (c1) at (\ox+\h*2,\oy,\oz+\w);
\coordinate (c2) at (\ox+\h*2,\oy+\w,\oz+\w);
\coordinate (c3) at (\ox+\h*2,\oy+\w,\oz);
\draw[dashed] (c1)--(c2)--(c3);


\def\ox{0}
\def\oy{0}
\def\oz{0}

\coordinate (b1) at (\ox,\oy,\oz+\h);
\coordinate (b2) at (\ox,\oy,\oz);
\coordinate (b3) at (\ox,\oy+\w,\oz);
\coordinate (b4) at (\ox,\oy+\w,\oz+\h);
\coordinate (b5) at (\ox+\w,\oy,\oz+\h);
\coordinate (b6) at (\ox+\w,\oy+\w,\oz+\h);
\coordinate (b7) at (\ox+\w,\oy+\w,\oz);
\draw[-,fill=white] (b1)--(b2)--(b3)--(b4)--(b1)--(b5)--(b6)--(b7)--(b3);
\draw[-] (b4)--(b6);
\coordinate (c1) at (\ox,\oy+\h,\oz);
\coordinate (c2) at (\ox,\oy+\h,\oz+\h);
\coordinate (c3) at (\ox+\w,\oy+\h,\oz+\h);
\draw[dashed] (c1)--(c2)--(c3);
\coordinate (c1) at (\ox,\oy+\h*2,\oz);
\coordinate (c2) at (\ox,\oy+\h*2,\oz+\h);
\coordinate (c3) at (\ox+\w,\oy+\h*2,\oz+\h);
\draw[dashed] (c1)--(c2)--(c3);
\coordinate (b1) at (\ox,\oy,\oz+\h+\os);
\coordinate (b2) at (\ox,\oy,\oz+\os);
\coordinate (b3) at (\ox,\oy+\w,\oz+\os);
\coordinate (b4) at (\ox,\oy+\w,\oz+\h+\os);
\coordinate (b5) at (\ox+\w,\oy,\oz+\h+\os);
\coordinate (b6) at (\ox+\w,\oy+\w,\oz+\h+\os);
\coordinate (b7) at (\ox+\w,\oy+\w,\oz+\os);
\draw[-,fill=white] (b1)--(b2)--(b3)--(b4)--(b1)--(b5)--(b6)--(b7)--(b3);
\draw[-] (b4)--(b6);
\draw[fill=white] (b1)--(b2)--(b3)--(b4)--(b1);
\coordinate (c1) at (\ox,\oy+\h,\oz+\os);
\coordinate (c2) at (\ox,\oy+\h,\oz+\h+\os);
\coordinate (c3) at (\ox+\w,\oy+\h,\oz+\h+\os);
\draw[dashed] (c1)--(c2)--(c3);
\coordinate (c1) at (\ox,\oy+\h*2,\oz+\os);
\coordinate (c2) at (\ox,\oy+\h*2,\oz+\h+\os);
\coordinate (c3) at (\ox+\w,\oy+\h*2,\oz+\h+\os);
\draw[dashed] (c1)--(c2)--(c3);
\coordinate (b1) at (\ox,\oy,\oz+\h+\os+\os);
\coordinate (b2) at (\ox,\oy,\oz+\os+\os);
\coordinate (b3) at (\ox,\oy+\w,\oz+\os+\os);
\coordinate (b4) at (\ox,\oy+\w,\oz+\h+\os+\os);
\coordinate (b5) at (\ox+\w,\oy,\oz+\h+\os+\os);
\coordinate (b6) at (\ox+\w,\oy+\w,\oz+\h+\os+\os);
\coordinate (b7) at (\ox+\w,\oy+\w,\oz+\os+\os);
\draw[-,fill=white] (b1)--(b2)--(b3)--(b4)--(b1)--(b5)--(b6)--(b7)--(b3);
\draw[-] (b4)--(b6);
\draw[fill=white] (b1)--(b2)--(b3)--(b4)--(b1);
\coordinate (c1) at (\ox,\oy+\h,\oz+\os+\os);
\coordinate (c2) at (\ox,\oy+\h,\oz+\h+\os+\os);
\coordinate (c3) at (\ox+\w,\oy+\h,\oz+\h+\os+\os);
\draw[dashed] (c1)--(c2)--(c3);
\coordinate (c1) at (\ox,\oy+\h*2,\oz+\os+\os);
\coordinate (c2) at (\ox,\oy+\h*2,\oz+\h+\os+\os);
\coordinate (c3) at (\ox+\w,\oy+\h*2,\oz+\h+\os+\os);
\draw[dashed] (c1)--(c2)--(c3);

 \node[circle,draw=black, fill=white, inner sep=0pt,minimum size=20pt] (node1) at (\ox+\ol,\oy+\ol,\oz+\ol) {{\Huge 1}};  


\def\ox{-5}
\def\oy{-5}
\def\oz{-7.5}

\coordinate (b1c1) at (\ox,\oy,\oz);
\coordinate (b1c2) at (\ox,\oy+\h,\oz);
\coordinate (b1c3) at (\ox,\oy+\h,\oz+\h);
\coordinate (b1c4) at (\ox,\oy,\oz+\h);
\draw[-] (b1c1)--(b1c2);
\draw[-] (b1c2)--(b1c3);
\draw[-] (b1c3)--(b1c4);
\draw[-] (b1c4)--(b1c1);
\coordinate (b1c1x) at (\ox+\h*3,\oy,\oz);
\coordinate (b1c2x) at (\ox+\h*3,\oy+\h,\oz);
\coordinate (b1c3x) at (\ox+\h*3,\oy+\h,\oz+\h);
\coordinate (b1c4x) at (\ox+\h*3,\oy,\oz+\h);
\draw[-] (b1c2x)--(b1c3x);
\draw[-] (b1c3x)--(b1c4x);
\draw [fill=white] (b1c3)--(b1c3x)--(b1c4x)--(b1c4)--(b1c3);
\draw [fill=white] (b1c3)--(b1c3x)--(b1c2x)--(b1c2)--(b1c3);

\coordinate (b2c1) at (\ox,\oy+\os,\oz);
\coordinate (b2c2) at (\ox,\oy+\h+\os,\oz);
\coordinate (b2c3) at (\ox,\oy+\h+\os,\oz+\h);
\coordinate (b2c4) at (\ox,\oy+\os,\oz+\h);
\draw[-] (b2c1)--(b2c2);
\draw[-] (b2c2)--(b2c3);
\draw[-] (b2c3)--(b2c4);
\draw[-] (b2c4)--(b2c1);
\coordinate (b2c1x) at (\ox+\h*3,\oy+\os,\oz);
\coordinate (b2c2x) at (\ox+\h*3,\oy+\h+\os,\oz);
\coordinate (b2c3x) at (\ox+\h*3,\oy+\h+\os,\oz+\h);
\coordinate (b2c4x) at (\ox+\h*3,\oy+\os,\oz+\h);
\draw[-] (b2c2x)--(b2c3x);
\draw[-] (b2c3x)--(b2c4x);
\draw [fill=white] (b2c3)--(b2c3x)--(b2c4x)--(b2c4)--(b2c3);
\draw [fill=white] (b2c3)--(b2c3x)--(b2c2x)--(b2c2)--(b2c3);
\draw [fill=white] (b2c1)--(b2c2)--(b2c3)--(b2c4)--(b2c1);

\coordinate (b3c1) at (\ox,\oy+\os+\os,\oz);
\coordinate (b3c2) at (\ox,\oy+\h+\os+\os,\oz);
\coordinate (b3c3) at (\ox,\oy+\h+\os+\os,\oz+\h);
\coordinate (b3c4) at (\ox,\oy+\os+\os,\oz+\h);
\draw[-] (b3c1)--(b3c2);
\draw[-] (b3c2)--(b3c3);
\draw[-] (b3c3)--(b3c4);
\draw[-] (b3c4)--(b3c1);
\coordinate (b3c1x) at (\ox+\h*3,\oy+\os+\os,\oz);
\coordinate (b3c2x) at (\ox+\h*3,\oy+\h+\os+\os,\oz);
\coordinate (b3c3x) at (\ox+\h*3,\oy+\h+\os+\os,\oz+\h);
\coordinate (b3c4x) at (\ox+\h*3,\oy+\os+\os,\oz+\h);
\draw[-] (b3c2x)--(b3c3x);
\draw[-] (b3c3x)--(b3c4x);
\draw [fill=white] (b3c3)--(b3c3x)--(b3c4x)--(b3c4)--(b3c3);
\draw [fill=white] (b3c3)--(b3c3x)--(b3c2x)--(b3c2)--(b3c3);
\draw [fill=white] (b3c1)--(b3c2)--(b3c3)--(b3c4)--(b3c1);

\coordinate (b4c1) at (\ox,\oy,\oz+\os);
\coordinate (b4c2) at (\ox,\oy+\h,\oz+\os);
\coordinate (b4c3) at (\ox,\oy+\h,\oz+\h+\os);
\coordinate (b4c4) at (\ox,\oy,\oz+\h+\os);
\draw[-] (b4c1)--(b4c2);
\draw[-] (b4c2)--(b4c3);
\draw[-] (b4c3)--(b4c4);
\draw[-] (b4c4)--(b4c1);
\coordinate (b4c1x) at (\ox+\h*3,\oy,\oz+\os);
\coordinate (b4c2x) at (\ox+\h*3,\oy+\h,\oz+\os);
\coordinate (b4c3x) at (\ox+\h*3,\oy+\h,\oz+\h+\os);
\coordinate (b4c4x) at (\ox+\h*3,\oy,\oz+\h+\os);
\draw[-] (b4c2x)--(b4c3x);
\draw[-] (b4c3x)--(b4c4x);
\draw [fill=white] (b4c3)--(b4c3x)--(b4c4x)--(b4c4)--(b4c3);
\draw [fill=white] (b4c3)--(b4c3x)--(b4c2x)--(b4c2)--(b4c3);
\draw [fill=white] (b4c1)--(b4c2)--(b4c3)--(b4c4)--(b4c1);

\coordinate (b5c1) at (\ox,\oy+\os,\oz+\os);
\coordinate (b5c2) at (\ox,\oy+\h+\os,\oz+\os);
\coordinate (b5c3) at (\ox,\oy+\h+\os,\oz+\h+\os);
\coordinate (b5c4) at (\ox,\oy+\os,\oz+\h+\os);
\draw[-] (b5c1)--(b5c2);
\draw[-] (b5c2)--(b5c3);
\draw[-] (b5c3)--(b5c4);
\draw[-] (b5c4)--(b5c1);
\coordinate (b5c1x) at (\ox+\h*3,\oy+\os,\oz+\os);
\coordinate (b5c2x) at (\ox+\h*3,\oy+\h+\os,\oz+\os);
\coordinate (b5c3x) at (\ox+\h*3,\oy+\h+\os,\oz+\h+\os);
\coordinate (b5c4x) at (\ox+\h*3,\oy+\os,\oz+\h+\os);
\draw[-] (b5c2x)--(b5c3x);
\draw[-] (b5c3x)--(b5c4x);
\draw [fill=white] (b5c3)--(b5c3x)--(b5c4x)--(b5c4)--(b5c3);
\draw [fill=white] (b5c3)--(b5c3x)--(b5c2x)--(b5c2)--(b5c3);
\draw [fill=white] (b5c1)--(b5c2)--(b5c3)--(b5c4)--(b5c1);

\coordinate (b6c1) at (\ox,\oy+\os+\os,\oz+\os);
\coordinate (b6c2) at (\ox,\oy+\h+\os+\os,\oz+\os);
\coordinate (b6c3) at (\ox,\oy+\h+\os+\os,\oz+\h+\os);
\coordinate (b6c4) at (\ox,\oy+\os+\os,\oz+\h+\os);
\draw[-] (b6c1)--(b6c2);
\draw[-] (b6c2)--(b6c3);
\draw[-] (b6c3)--(b6c4);
\draw[-] (b6c4)--(b6c1);
\coordinate (b6c1x) at (\ox+\h*3,\oy+\os+\os,\oz+\os);
\coordinate (b6c2x) at (\ox+\h*3,\oy+\h+\os+\os,\oz+\os);
\coordinate (b6c3x) at (\ox+\h*3,\oy+\h+\os+\os,\oz+\h+\os);
\coordinate (b6c4x) at (\ox+\h*3,\oy+\os+\os,\oz+\h+\os);
\draw[-] (b6c2x)--(b6c3x);
\draw[-] (b6c3x)--(b6c4x);
\draw [fill=white] (b6c3)--(b6c3x)--(b6c4x)--(b6c4)--(b6c3);
\draw [fill=white] (b6c3)--(b6c3x)--(b6c2x)--(b6c2)--(b6c3);
\draw [fill=white] (b6c1)--(b6c2)--(b6c3)--(b6c4)--(b6c1);

\coordinate (b7c1) at (\ox,\oy,\oz+\os+\os);
\coordinate (b7c2) at (\ox,\oy+\h,\oz+\os+\os);
\coordinate (b7c3) at (\ox,\oy+\h,\oz+\h+\os+\os);
\coordinate (b7c4) at (\ox,\oy,\oz+\h+\os+\os);
\draw[-] (b7c1)--(b7c2);
\draw[-] (b7c2)--(b7c3);
\draw[-] (b7c3)--(b7c4);
\draw[-] (b7c4)--(b7c1);
\coordinate (b7c1x) at (\ox+\h*3,\oy,\oz+\os+\os);
\coordinate (b7c2x) at (\ox+\h*3,\oy+\h,\oz+\os+\os);
\coordinate (b7c3x) at (\ox+\h*3,\oy+\h,\oz+\h+\os+\os);
\coordinate (b7c4x) at (\ox+\h*3,\oy,\oz+\h+\os+\os);
\draw[-] (b7c2x)--(b7c3x);
\draw[-] (b7c3x)--(b7c4x);
\draw [fill=white] (b7c3)--(b7c3x)--(b7c4x)--(b7c4)--(b7c3);
\draw [fill=white] (b7c3)--(b7c3x)--(b7c2x)--(b7c2)--(b7c3);
\draw [fill=white] (b7c1)--(b7c2)--(b7c3)--(b7c4)--(b7c1);

\coordinate (b8c1) at (\ox,\oy+\os,\oz+\os+\os);
\coordinate (b8c2) at (\ox,\oy+\h+\os,\oz+\os+\os);
\coordinate (b8c3) at (\ox,\oy+\h+\os,\oz+\h+\os+\os);
\coordinate (b8c4) at (\ox,\oy+\os,\oz+\h+\os+\os);
\draw[-] (b8c1)--(b8c2);
\draw[-] (b8c2)--(b8c3);
\draw[-] (b8c3)--(b8c4);
\draw[-] (b8c4)--(b8c1);
\coordinate (b8c1x) at (\ox+\h*3,\oy+\os,\oz+\os+\os);
\coordinate (b8c2x) at (\ox+\h*3,\oy+\h+\os,\oz+\os+\os);
\coordinate (b8c3x) at (\ox+\h*3,\oy+\h+\os,\oz+\h+\os+\os);
\coordinate (b8c4x) at (\ox+\h*3,\oy+\os,\oz+\h+\os+\os);
\draw[-] (b8c2x)--(b8c3x);
\draw[-] (b8c3x)--(b8c4x);
\draw [fill=white] (b8c3)--(b8c3x)--(b8c4x)--(b8c4)--(b8c3);
\draw [fill=white] (b8c3)--(b8c3x)--(b8c2x)--(b8c2)--(b8c3);
\draw [fill=white] (b8c1)--(b8c2)--(b8c3)--(b8c4)--(b8c1);

\coordinate (b9c1) at (\ox,\oy+\os+\os,\oz+\os+\os);
\coordinate (b9c2) at (\ox,\oy+\h+\os+\os,\oz+\os+\os);
\coordinate (b9c3) at (\ox,\oy+\h+\os+\os,\oz+\h+\os+\os);
\coordinate (b9c4) at (\ox,\oy+\os+\os,\oz+\h+\os+\os);
\draw[-] (b9c1)--(b9c2);
\draw[-] (b9c2)--(b9c3);
\draw[-] (b9c3)--(b9c4);
\draw[-] (b9c4)--(b9c1);
\coordinate (b9c1x) at (\ox+\h*3,\oy+\os+\os,\oz+\os+\os);
\coordinate (b9c2x) at (\ox+\h*3,\oy+\h+\os+\os,\oz+\os+\os);
\coordinate (b9c3x) at (\ox+\h*3,\oy+\h+\os+\os,\oz+\h+\os+\os);
\coordinate (b9c4x) at (\ox+\h*3,\oy+\os+\os,\oz+\h+\os+\os);
\draw[-] (b9c2x)--(b9c3x);
\draw[-] (b9c3x)--(b9c4x);
\draw [fill=white] (b9c3)--(b9c3x)--(b9c4x)--(b9c4)--(b9c3);
\draw [fill=white] (b9c3)--(b9c3x)--(b9c2x)--(b9c2)--(b9c3);
\draw [fill=white] (b9c1)--(b9c2)--(b9c3)--(b9c4)--(b9c1);

 \node[circle,draw=black, fill=white, inner sep=0pt,minimum size=20pt] (node2) at (\ox+\ol,\oy+\ol,\oz+\ol) {{\Huge 2}};  



\def\ox{-.25}
\def\oy{-.25}
\def\oz{-7.5}

\coordinate (b1) at (\ox,\oy,\oz+\w);
\coordinate (b2) at (\ox,\oy,\oz);
\coordinate (b3) at (\ox,\oy+\h,\oz);
\coordinate (b4) at (\ox,\oy+\h,\oz+\w);
\coordinate (b5) at (\ox+\w,\oy,\oz+\w);
\coordinate (b6) at (\ox+\w,\oy+\h,\oz+\w);
\coordinate (b7) at (\ox+\w,\oy+\h,\oz);
\draw[-,fill=white] (b1)--(b2)--(b3)--(b4)--(b1)--(b5)--(b6)--(b7)--(b3);
\draw[-] (b4)--(b6);

\coordinate (c1) at (\ox+\h,\oy,\oz+\w);
\coordinate (c2) at (\ox+\h,\oy+\h,\oz+\w);
\coordinate (c3) at (\ox+\h,\oy+\h,\oz);
\draw[dashed] (c1)--(c2)--(c3);

\coordinate (c1) at (\ox+\h*2,\oy,\oz+\w);
\coordinate (c2) at (\ox+\h*2,\oy+\h,\oz+\w);
\coordinate (c3) at (\ox+\h*2,\oy+\h,\oz);
\draw[dashed] (c1)--(c2)--(c3);

\coordinate (b1) at (\ox,\oy+\os,\oz+\w);
\coordinate (b2) at (\ox,\oy+\os,\oz);
\coordinate (b3) at (\ox,\oy+\h+\os,\oz);
\coordinate (b4) at (\ox,\oy+\h+\os,\oz+\w);
\coordinate (b5) at (\ox+\w,\oy+\os,\oz+\w);
\coordinate (b6) at (\ox+\w,\oy+\h+\os,\oz+\w);
\coordinate (b7) at (\ox+\w,\oy+\h+\os,\oz);
\draw[fill=white] (b1)--(b2)--(b3)--(b7)--(b6)--(b5)--(b1);
\draw[-] (b6)--(b4)--(b1);
\draw[-] (b4)--(b3);

\coordinate (c1) at (\ox+\h,\oy+\os,\oz+\w);
\coordinate (c2) at (\ox+\h,\oy+\h+\os,\oz+\w);
\coordinate (c3) at (\ox+\h,\oy+\h+\os,\oz);
\draw[dashed] (c1)--(c2)--(c3);

\coordinate (c1) at (\ox+\h*2,\oy+\os,\oz+\w);
\coordinate (c2) at (\ox+\h*2,\oy+\h+\os,\oz+\w);
\coordinate (c3) at (\ox+\h*2,\oy+\h+\os,\oz);
\draw[dashed] (c1)--(c2)--(c3);

\coordinate (b1) at (\ox,\oy+\os+\os,\oz+\w);
\coordinate (b2) at (\ox,\oy+\os+\os,\oz);
\coordinate (b3) at (\ox,\oy+\h+\os+\os,\oz);
\coordinate (b4) at (\ox,\oy+\h+\os+\os,\oz+\w);
\coordinate (b5) at (\ox+\w,\oy+\os+\os,\oz+\w);
\coordinate (b6) at (\ox+\w,\oy+\h+\os+\os,\oz+\w);
\coordinate (b7) at (\ox+\w,\oy+\h+\os+\os,\oz);
\draw[fill=white] (b1)--(b2)--(b3)--(b7)--(b6)--(b5)--(b1);
\draw[-] (b6)--(b4)--(b1);
\draw[-] (b4)--(b3);

\coordinate (c1) at (\ox+\h,\oy+\os+\os,\oz+\w);
\coordinate (c2) at (\ox+\h,\oy+\h+\os+\os,\oz+\w);
\coordinate (c3) at (\ox+\h,\oy+\h+\os+\os,\oz);
\draw[dashed] (c1)--(c2)--(c3);

\coordinate (c1) at (\ox+\h*2,\oy+\os+\os,\oz+\w);
\coordinate (c2) at (\ox+\h*2,\oy+\h+\os+\os,\oz+\w);
\coordinate (c3) at (\ox+\h*2,\oy+\h+\os+\os,\oz);
\draw[dashed] (c1)--(c2)--(c3);

 \node[circle,draw=black, fill=white, inner sep=0pt,minimum size=20pt] () at (\ox+\ol,\oy+\ol,\oz+\ol) {{\Huge 3}};  



\def\ox{-5}
\def\oy{-5}
\def\oz{-15}

\coordinate (b1c1) at (\ox,\oy,\oz);
\coordinate (b1c2) at (\ox,\oy+\h,\oz);
\coordinate (b1c3) at (\ox,\oy+\h,\oz+\h);
\coordinate (b1c4) at (\ox,\oy,\oz+\h);
\draw[-] (b1c1)--(b1c2);
\draw[-] (b1c2)--(b1c3);
\draw[-] (b1c3)--(b1c4);
\draw[-] (b1c4)--(b1c1);
\coordinate (b1c1x) at (\ox+\h*3,\oy,\oz);
\coordinate (b1c2x) at (\ox+\h*3,\oy+\h,\oz);
\coordinate (b1c3x) at (\ox+\h*3,\oy+\h,\oz+\h);
\coordinate (b1c4x) at (\ox+\h*3,\oy,\oz+\h);
\draw[-] (b1c2x)--(b1c3x);
\draw[-] (b1c3x)--(b1c4x);
\draw [fill=white] (b1c3)--(b1c3x)--(b1c4x)--(b1c4)--(b1c3);
\draw [fill=white] (b1c3)--(b1c3x)--(b1c2x)--(b1c2)--(b1c3);

\coordinate (b2c1) at (\ox,\oy+\os,\oz);
\coordinate (b2c2) at (\ox,\oy+\h+\os,\oz);
\coordinate (b2c3) at (\ox,\oy+\h+\os,\oz+\h);
\coordinate (b2c4) at (\ox,\oy+\os,\oz+\h);
\draw[-] (b2c1)--(b2c2);
\draw[-] (b2c2)--(b2c3);
\draw[-] (b2c3)--(b2c4);
\draw[-] (b2c4)--(b2c1);
\coordinate (b2c1x) at (\ox+\h*3,\oy+\os,\oz);
\coordinate (b2c2x) at (\ox+\h*3,\oy+\h+\os,\oz);
\coordinate (b2c3x) at (\ox+\h*3,\oy+\h+\os,\oz+\h);
\coordinate (b2c4x) at (\ox+\h*3,\oy+\os,\oz+\h);
\draw[-] (b2c2x)--(b2c3x);
\draw[-] (b2c3x)--(b2c4x);
\draw [fill=white] (b2c3)--(b2c3x)--(b2c4x)--(b2c4)--(b2c3);
\draw [fill=white] (b2c3)--(b2c3x)--(b2c2x)--(b2c2)--(b2c3);
\draw [fill=white] (b2c1)--(b2c2)--(b2c3)--(b2c4)--(b2c1);

\coordinate (b3c1) at (\ox,\oy+\os+\os,\oz);
\coordinate (b3c2) at (\ox,\oy+\h+\os+\os,\oz);
\coordinate (b3c3) at (\ox,\oy+\h+\os+\os,\oz+\h);
\coordinate (b3c4) at (\ox,\oy+\os+\os,\oz+\h);
\draw[-] (b3c1)--(b3c2);
\draw[-] (b3c2)--(b3c3);
\draw[-] (b3c3)--(b3c4);
\draw[-] (b3c4)--(b3c1);
\coordinate (b3c1x) at (\ox+\h*3,\oy+\os+\os,\oz);
\coordinate (b3c2x) at (\ox+\h*3,\oy+\h+\os+\os,\oz);
\coordinate (b3c3x) at (\ox+\h*3,\oy+\h+\os+\os,\oz+\h);
\coordinate (b3c4x) at (\ox+\h*3,\oy+\os+\os,\oz+\h);
\draw[-] (b3c2x)--(b3c3x);
\draw[-] (b3c3x)--(b3c4x);
\draw [fill=white] (b3c3)--(b3c3x)--(b3c4x)--(b3c4)--(b3c3);
\draw [fill=white] (b3c3)--(b3c3x)--(b3c2x)--(b3c2)--(b3c3);
\draw [fill=white] (b3c1)--(b3c2)--(b3c3)--(b3c4)--(b3c1);

\coordinate (b4c1) at (\ox,\oy,\oz+\os);
\coordinate (b4c2) at (\ox,\oy+\h,\oz+\os);
\coordinate (b4c3) at (\ox,\oy+\h,\oz+\h+\os);
\coordinate (b4c4) at (\ox,\oy,\oz+\h+\os);
\draw[-] (b4c1)--(b4c2);
\draw[-] (b4c2)--(b4c3);
\draw[-] (b4c3)--(b4c4);
\draw[-] (b4c4)--(b4c1);
\coordinate (b4c1x) at (\ox+\h*3,\oy,\oz+\os);
\coordinate (b4c2x) at (\ox+\h*3,\oy+\h,\oz+\os);
\coordinate (b4c3x) at (\ox+\h*3,\oy+\h,\oz+\h+\os);
\coordinate (b4c4x) at (\ox+\h*3,\oy,\oz+\h+\os);
\draw[-] (b4c2x)--(b4c3x);
\draw[-] (b4c3x)--(b4c4x);
\draw [fill=white] (b4c3)--(b4c3x)--(b4c4x)--(b4c4)--(b4c3);
\draw [fill=white] (b4c3)--(b4c3x)--(b4c2x)--(b4c2)--(b4c3);
\draw [fill=white] (b4c1)--(b4c2)--(b4c3)--(b4c4)--(b4c1);

\coordinate (b5c1) at (\ox,\oy+\os,\oz+\os);
\coordinate (b5c2) at (\ox,\oy+\h+\os,\oz+\os);
\coordinate (b5c3) at (\ox,\oy+\h+\os,\oz+\h+\os);
\coordinate (b5c4) at (\ox,\oy+\os,\oz+\h+\os);
\draw[-] (b5c1)--(b5c2);
\draw[-] (b5c2)--(b5c3);
\draw[-] (b5c3)--(b5c4);
\draw[-] (b5c4)--(b5c1);
\coordinate (b5c1x) at (\ox+\h*3,\oy+\os,\oz+\os);
\coordinate (b5c2x) at (\ox+\h*3,\oy+\h+\os,\oz+\os);
\coordinate (b5c3x) at (\ox+\h*3,\oy+\h+\os,\oz+\h+\os);
\coordinate (b5c4x) at (\ox+\h*3,\oy+\os,\oz+\h+\os);
\draw[-] (b5c2x)--(b5c3x);
\draw[-] (b5c3x)--(b5c4x);
\draw [fill=white] (b5c3)--(b5c3x)--(b5c4x)--(b5c4)--(b5c3);
\draw [fill=white] (b5c3)--(b5c3x)--(b5c2x)--(b5c2)--(b5c3);
\draw [fill=white] (b5c1)--(b5c2)--(b5c3)--(b5c4)--(b5c1);

\coordinate (b6c1) at (\ox,\oy+\os+\os,\oz+\os);
\coordinate (b6c2) at (\ox,\oy+\h+\os+\os,\oz+\os);
\coordinate (b6c3) at (\ox,\oy+\h+\os+\os,\oz+\h+\os);
\coordinate (b6c4) at (\ox,\oy+\os+\os,\oz+\h+\os);
\draw[-] (b6c1)--(b6c2);
\draw[-] (b6c2)--(b6c3);
\draw[-] (b6c3)--(b6c4);
\draw[-] (b6c4)--(b6c1);
\coordinate (b6c1x) at (\ox+\h*3,\oy+\os+\os,\oz+\os);
\coordinate (b6c2x) at (\ox+\h*3,\oy+\h+\os+\os,\oz+\os);
\coordinate (b6c3x) at (\ox+\h*3,\oy+\h+\os+\os,\oz+\h+\os);
\coordinate (b6c4x) at (\ox+\h*3,\oy+\os+\os,\oz+\h+\os);
\draw[-] (b6c2x)--(b6c3x);
\draw[-] (b6c3x)--(b6c4x);
\draw [fill=white] (b6c3)--(b6c3x)--(b6c4x)--(b6c4)--(b6c3);
\draw [fill=white] (b6c3)--(b6c3x)--(b6c2x)--(b6c2)--(b6c3);
\draw [fill=white] (b6c1)--(b6c2)--(b6c3)--(b6c4)--(b6c1);

\coordinate (b7c1) at (\ox,\oy,\oz+\os+\os);
\coordinate (b7c2) at (\ox,\oy+\h,\oz+\os+\os);
\coordinate (b7c3) at (\ox,\oy+\h,\oz+\h+\os+\os);
\coordinate (b7c4) at (\ox,\oy,\oz+\h+\os+\os);
\draw[-] (b7c1)--(b7c2);
\draw[-] (b7c2)--(b7c3);
\draw[-] (b7c3)--(b7c4);
\draw[-] (b7c4)--(b7c1);
\coordinate (b7c1x) at (\ox+\h*3,\oy,\oz+\os+\os);
\coordinate (b7c2x) at (\ox+\h*3,\oy+\h,\oz+\os+\os);
\coordinate (b7c3x) at (\ox+\h*3,\oy+\h,\oz+\h+\os+\os);
\coordinate (b7c4x) at (\ox+\h*3,\oy,\oz+\h+\os+\os);
\draw[-] (b7c2x)--(b7c3x);
\draw[-] (b7c3x)--(b7c4x);
\draw [fill=white] (b7c3)--(b7c3x)--(b7c4x)--(b7c4)--(b7c3);
\draw [fill=white] (b7c3)--(b7c3x)--(b7c2x)--(b7c2)--(b7c3);
\draw [fill=white] (b7c1)--(b7c2)--(b7c3)--(b7c4)--(b7c1);

\coordinate (b8c1) at (\ox,\oy+\os,\oz+\os+\os);
\coordinate (b8c2) at (\ox,\oy+\h+\os,\oz+\os+\os);
\coordinate (b8c3) at (\ox,\oy+\h+\os,\oz+\h+\os+\os);
\coordinate (b8c4) at (\ox,\oy+\os,\oz+\h+\os+\os);
\draw[-] (b8c1)--(b8c2);
\draw[-] (b8c2)--(b8c3);
\draw[-] (b8c3)--(b8c4);
\draw[-] (b8c4)--(b8c1);
\coordinate (b8c1x) at (\ox+\h*3,\oy+\os,\oz+\os+\os);
\coordinate (b8c2x) at (\ox+\h*3,\oy+\h+\os,\oz+\os+\os);
\coordinate (b8c3x) at (\ox+\h*3,\oy+\h+\os,\oz+\h+\os+\os);
\coordinate (b8c4x) at (\ox+\h*3,\oy+\os,\oz+\h+\os+\os);
\draw[-] (b8c2x)--(b8c3x);
\draw[-] (b8c3x)--(b8c4x);
\draw [fill=white] (b8c3)--(b8c3x)--(b8c4x)--(b8c4)--(b8c3);
\draw [fill=white] (b8c3)--(b8c3x)--(b8c2x)--(b8c2)--(b8c3);
\draw [fill=white] (b8c1)--(b8c2)--(b8c3)--(b8c4)--(b8c1);

\coordinate (b9c1) at (\ox,\oy+\os+\os,\oz+\os+\os);
\coordinate (b9c2) at (\ox,\oy+\h+\os+\os,\oz+\os+\os);
\coordinate (b9c3) at (\ox,\oy+\h+\os+\os,\oz+\h+\os+\os);
\coordinate (b9c4) at (\ox,\oy+\os+\os,\oz+\h+\os+\os);
\draw[-] (b9c1)--(b9c2);
\draw[-] (b9c2)--(b9c3);
\draw[-] (b9c3)--(b9c4);
\draw[-] (b9c4)--(b9c1);
\coordinate (b9c1x) at (\ox+\h*3,\oy+\os+\os,\oz+\os+\os);
\coordinate (b9c2x) at (\ox+\h*3,\oy+\h+\os+\os,\oz+\os+\os);
\coordinate (b9c3x) at (\ox+\h*3,\oy+\h+\os+\os,\oz+\h+\os+\os);
\coordinate (b9c4x) at (\ox+\h*3,\oy+\os+\os,\oz+\h+\os+\os);
\draw[-] (b9c2x)--(b9c3x);
\draw[-] (b9c3x)--(b9c4x);
\draw [fill=white] (b9c3)--(b9c3x)--(b9c4x)--(b9c4)--(b9c3);
\draw [fill=white] (b9c3)--(b9c3x)--(b9c2x)--(b9c2)--(b9c3);
\draw [fill=white] (b9c1)--(b9c2)--(b9c3)--(b9c4)--(b9c1);

 \node[circle,draw=black, fill=white, inner sep=0pt,minimum size=20pt] () at (\ox+\ol,\oy+\ol,\oz+\ol) {{\Huge 4}};  



\def\ox{0}
\def\oy{0}
\def\oz{-15}

\coordinate (b1) at (\ox,\oy,\oz+\h);
\coordinate (b2) at (\ox,\oy,\oz);
\coordinate (b3) at (\ox,\oy+\w,\oz);
\coordinate (b4) at (\ox,\oy+\w,\oz+\h);
\coordinate (b5) at (\ox+\w,\oy,\oz+\h);
\coordinate (b6) at (\ox+\w,\oy+\w,\oz+\h);
\coordinate (b7) at (\ox+\w,\oy+\w,\oz);
\draw[-,fill=white] (b1)--(b2)--(b3)--(b4)--(b1)--(b5)--(b6)--(b7)--(b3);
\draw[-] (b4)--(b6);
\coordinate (c1) at (\ox,\oy+\h,\oz);
\coordinate (c2) at (\ox,\oy+\h,\oz+\h);
\coordinate (c3) at (\ox+\w,\oy+\h,\oz+\h);
\draw[dashed] (c1)--(c2)--(c3);
\coordinate (c1) at (\ox,\oy+\h*2,\oz);
\coordinate (c2) at (\ox,\oy+\h*2,\oz+\h);
\coordinate (c3) at (\ox+\w,\oy+\h*2,\oz+\h);
\draw[dashed] (c1)--(c2)--(c3);
\coordinate (b1) at (\ox,\oy,\oz+\h+\os);
\coordinate (b2) at (\ox,\oy,\oz+\os);
\coordinate (b3) at (\ox,\oy+\w,\oz+\os);
\coordinate (b4) at (\ox,\oy+\w,\oz+\h+\os);
\coordinate (b5) at (\ox+\w,\oy,\oz+\h+\os);
\coordinate (b6) at (\ox+\w,\oy+\w,\oz+\h+\os);
\coordinate (b7) at (\ox+\w,\oy+\w,\oz+\os);
\draw[-,fill=white] (b1)--(b2)--(b3)--(b4)--(b1)--(b5)--(b6)--(b7)--(b3);
\draw[-] (b4)--(b6);
\draw[fill=white] (b1)--(b2)--(b3)--(b4)--(b1);
\coordinate (c1) at (\ox,\oy+\h,\oz+\os);
\coordinate (c2) at (\ox,\oy+\h,\oz+\h+\os);
\coordinate (c3) at (\ox+\w,\oy+\h,\oz+\h+\os);
\draw[dashed] (c1)--(c2)--(c3);
\coordinate (c1) at (\ox,\oy+\h*2,\oz+\os);
\coordinate (c2) at (\ox,\oy+\h*2,\oz+\h+\os);
\coordinate (c3) at (\ox+\w,\oy+\h*2,\oz+\h+\os);
\draw[dashed] (c1)--(c2)--(c3);
\coordinate (b1) at (\ox,\oy,\oz+\h+\os+\os);
\coordinate (b2) at (\ox,\oy,\oz+\os+\os);
\coordinate (b3) at (\ox,\oy+\w,\oz+\os+\os);
\coordinate (b4) at (\ox,\oy+\w,\oz+\h+\os+\os);
\coordinate (b5) at (\ox+\w,\oy,\oz+\h+\os+\os);
\coordinate (b6) at (\ox+\w,\oy+\w,\oz+\h+\os+\os);
\coordinate (b7) at (\ox+\w,\oy+\w,\oz+\os+\os);
\draw[-,fill=white] (b1)--(b2)--(b3)--(b4)--(b1)--(b5)--(b6)--(b7)--(b3);
\draw[-] (b4)--(b6);
\draw[fill=white] (b1)--(b2)--(b3)--(b4)--(b1);
\coordinate (c1) at (\ox,\oy+\h,\oz+\os+\os);
\coordinate (c2) at (\ox,\oy+\h,\oz+\h+\os+\os);
\coordinate (c3) at (\ox+\w,\oy+\h,\oz+\h+\os+\os);
\draw[dashed] (c1)--(c2)--(c3);
\coordinate (c1) at (\ox,\oy+\h*2,\oz+\os+\os);
\coordinate (c2) at (\ox,\oy+\h*2,\oz+\h+\os+\os);
\coordinate (c3) at (\ox+\w,\oy+\h*2,\oz+\h+\os+\os);
\draw[dashed] (c1)--(c2)--(c3);

\node[circle,draw=black, fill=white, inner sep=0pt,minimum size=20pt] () at (\ox+\ol,\oy+\ol,\oz+\ol) {{\Huge 5}};



\end{tikzpicture}
\end{center}

\label{mpi_send_image}
\end{figure}

The communication between the MPI processes is presented in Figure \ref{mpi_send_image}. The key details of the implementation can be considered as follows. If there are $np$ MPI processes available for parallel execution, then on the first step of the parallel algorithm, each process performs $kpz = \lfloor Nz/np \rfloor$ or $ kpz = \lfloor Nz/np \rfloor + 1$ 2D DSTs of the 2D horizontal slices of the 3D array of the right-hand side. After completion of the $1st$ step, each process will send a part of the 3D array of the transformed data to every other process. The size of the submitted data to another process is 
 $ N_x \times kpy \times kpz $, where $ kpy = \lfloor N_y/np \rfloor $ or $ kpy = \lfloor N_y/np \rfloor + 1 $. The second step of the process is the solution of the $N_x \times N_y$ independent tridiagonal linear systems of the size $N_z \times N_z$. In this case, the 3D array of the transformed right-hand side is divided along the $y-$ direction and every process has to solve $ N_x \times kpy$ systems.   As a result of the second step, the program obtains the transformed solution of the system. Then an individual process sends a portion of this array, of the size $N_x \times kpy \times kpz $, to every other process to set up the last inverse transformation step. This step is executed in the same way as the forward transformation except for the use of the transformed solution array rather than the right-hand side.

Despite its immense capability, the MPI implementation's performance is limited in this approach. Since the algorithm uses 2D DST, the program distributes sets of the 2D slices of the 3D working arrays between the processes. This significantly reduces the communication time between the processes when compared to the transformation step that uses 
parallelization of the sets of 1D FFTs.


However, this approach distributes the 2D DST operations in the $z-$ direction of the computational domain in both the forward and inverse transforms. These are the most time expensive steps of the solver. This means that the number of processes that can be used in parallel is bounded by $N_z$, the number of grid points in the $z-$ direction of the computational domain.  

\subsection{Hybrid} \label{sec: Hybrid}

This subsection discusses an approach that combines the advantages of both OpenMP and MPI tools in the presented algorithm. Consider a cluster with $NH$ nodes each with $KH$ cores. An OpenMP program can only be run on a single node, so only $KH$ threads can be used for the parallel calculations. As previously mentioned, the number of MPI processes is bounded by $Nz$ in the MPI implementation. To use the full power of a cluster, that is to utilize all $NH \times KH$ available cores, it is possible to combine both the OpenMP and MPI tools into a hybrid program. MPI non-blocking $Isend$ and $Irecv$ commands can be used to transfer data between the nodes where the shared memory is used by OpenMP threads to implement allocated tasks in parallel. This has a clear advantage over using the strictly MPI approach as it will reduce time lost to communication between MPI processes.

In this hybrid case, an MPI process uses an entire node. Then OpenMP allows access to every core on the node. As in the case of the MPI implementation, the computational domain needs to be divided along the $z-$ direction, this number is still bounded by $N_z$. However, the hybrid approach uses all available cores on a node via OpenMP threads to parallelize the 2D DST. This is accomplished by using FFTW multi-threading, see \cite{FFTW_doc}. The implementation of this modified approach is outlined in Algorithm \ref{Hybrid}.

\begin{algorithm}[H]
\caption{Hybrid 3D Helmholtz Direct Solver}
\begin{algorithmic} [1] \label{Hybrid}
\STATE Create multi-threaded FFTW plan
\STATE Find $start_y,start_z,end_y,$ and $end_z$ using the rank
\FOR{$l=start_z,\dots,end_z$}
\STATE 2D forward multi-threading DST in $x-,y-$ direction
\ENDFOR
\STATE Scatter the data to the appropriate process
\STATE \#pragma omp parallel for collapse(2)
\FOR{$j=start_y,\dots,end_y;i=1,\dots,N_x$}
\STATE Solve the tridiagonal system using LU decomposition
\ENDFOR
\STATE Scatter the data to the appropriate process	
\FOR{$l=start_z,\dots,end_z$}
\STATE 2D inverse multi-threading DST in $x-,y-$ direction
\ENDFOR
\end{algorithmic}
\end{algorithm}

\section{Numerical Results} \label{results_section}

In this section, the results of numerical experiments that demonstrate the quality of the proposed numerical methods are presented. These algorithms were implemented in C programming language and the majority of the numerical experiments were conducted on the ``Cori'' cluster at Lawrence Berkeley National Laboratory with Haswell nodes. The Haswell nodes contain 32 Intel Xeon Haswell processors with approximately 2.3 GHz clock frequency. For comparison with previously published results, we also considered several experiments on a standard iMac desktop with an Intel Core i7, 2.93 GHz processor and 16 Gb of RAM, and a Xeon X5690 server running at 3.47 GHz with 144 Gb of RAM.

\subsection{Sequential Implementation of the Direct FFT Solvers}

First, we investigate the efficiency of the developed direct solvers in the case of a 3D test problem and sequential implementation. We choose to illustrate the quality of the developed direct methods on the test problems recently published in \cite{Turkel}. The authors of the paper considered the solution of several 3D test problems using iterative block-parallel CARP-CG method \citep{Gordon}. We calculate the solution of the same problems by applying direct algorithms discussed in the previous section. Also, we demonstrate that even the sequential variant of the developed method is significantly faster than the mentioned iterative solver when implemented on less expensive hardware. But this could be expected since the iterative CARP-CG method is designed for the solution of the general 3D Helmholtz equation instead of the problems with specific restrictions on the problem coefficient considered in this paper.

In these test problems, the coefficient $k(z) = a - b\sin(cz)$ with $ a > b \ge 0$ depends only on one spatial variable, i.e.\,the developed methods could be used as direct solvers to find an approximate solution to the boundary value problem (\ref{problem}, \ref{bc}). In our experiments, we use the following measures related to the approximate $U$ and analytic $u$ solutions of the problems:
\begin{list}{$\bullet$}{}
\item \textbf{$L_2$-res} is $\|AU-F\|_{2}$,   
\item \textbf{$L_2$-err}(the relative $L_2$ error) is $\|u-U\|_2/\|u\|_2$.   
\item \textbf{max-err} is $\|u-U\|_{\infty}$.   
\end{list}

\noindent In the following numerical experiments the analytic solution 
\begin{equation}
u(x,y,z)=\sin(\beta x)\sin(\gamma y) e^{-\frac{k(z)}{c} } , \textnormal {where }  \beta^2 + \gamma^2 = a^2+b^2  \nonumber  \\
\end{equation}
is used. We also assume that $ L_{\alpha}^l = 0 $ and $ L_{\alpha}^r=\pi, \alpha=x,y,z. $ This is the same solution considered in \cite{Turkel} up to the notation for the independent variables. Then the right-hand side of (\ref{problem}) is $ f(x,y,z)= -b(2a+c)\sin(cz)e^{-\frac{k(z)}{c}}\sin(\beta x)\sin(\gamma y).$ We consider the application of the developed direct methods to the solutions obtained by the iterative approach used in \cite{Turkel} on the Supermicro cluster consisting of 12 nodes. Each node had two Intel Xeon E5520 quad CPUs running at 2.27 GHz. The two CPUs shared 8 GB of memory. We restrict our consideration to the numerical results with available CPU time and corresponding to the smallest presented in  \cite{Turkel} relative $L_2$ error of 0.001. One can find these results in Tables 2 and 3 in the mentioned paper. For a demonstration of the efficiency of the presented direct solvers, we run all test problems on both standard desktop iMac i7, and on the Xeon X5690 server.

In the first experiment, we use the following parameters: $a=10, b=9, \gamma = 9 \ \ ( 1 \le  k \le 19)$. Table 1 presents a comparison of various solvers:
the first two rows show the iterative solver used in \cite{Turkel}. Rows from 3 to 8 present results of the second-order direct solver considered in our previous publications \cite{gkl,yg2}, and fourth and sixth-order solvers presented in Section 3 of this paper. Rows 3-5 give results for the Xeon X5690 server and lines 6-8 exhibit the results achieved on an iMac PC. The first column represents the hardware used in the numerical experiment. The second and third columns indicate the order of approximation of the solver and the type of the solver (direct or iterative). In the fourth column, the number of grid points needed to reach the indicated relative accuracy (\textbf{$L_2$-err} $< 0.001$). The fifth column shows the number of iterations until the convergence of the iterative solver, where the case of the direct solvers we put 1. The last column displays the CPU time required for each test run.       

\begin{table}[H]
\centering
\captionsetup{justification=centering}
\begin{tabular}{|r | r | r | r | r | r | }
\hline
     CPU  & Scheme & Type &$N$ & $ \#$ iter.  & Time(s) \\  \hline
SM cl & 2& iter.& 333 &1970& 703 \\ \hline
SM cl & 6& iter.&45 &350& 1.01  \\ \hline
X5690 & 2&dir.& 353 &1&15.18   \\ \hline
X5690 & 4&dir.& 62& 1 & 0.078  \\ \hline
X5690 & 6&dir.& 50 & 1& .055   \\ \hline
iMac i7 & 2&dir.& 353 &1& 19.8  \\ \hline
iMac i7 & 4&dir.&  62 &1& .097  \\ \hline
iMac i7 & 6&dir.& 50 &1& .08  \\ \hline
\end{tabular}
\caption{Comparison of Direct and Iterative Solvers on the First 3D Test Problem}
\end{table}

In our experiments, to reach the desired accuracy with the second-order approximation scheme, we needed to use a $353^3$ grid, and the direct solver on X5690 and iMac i7 were 46 and 36 times faster than the iterative solver, respectively. In the case of the sixth-order scheme, the direct solver on X5690 and iMac i7 were, 18 and 13 times faster than the iterative solver, respectively. We must mention that the CPU time for the fourth-order scheme on $64^3$ grid was $0.07 sec$ on X5690. These numbers indicate that due to the optimality condition of the FFT method, sometimes it is advantageous to consider a slightly larger number of grid points which has more factors of 2 in its prime factorization.

In the next experiment, we consider the same problem with parameters  $a=80, b=40, \gamma = 40 \ \ ( 40 \le  k \le120)$ and various values for $c =10, 50, 70, 80$. As in the previous series of numerical experiments, we only consider the solutions of the test problems to reach \textbf{$L_2$-err} = 0.001. In \cite{Turkel}, there is no data for the second-order scheme since it was stated that ``the second-order scheme could not achieve the error goals with grids of manageable sizes'' We also restrict our consideration to the sixth-order direct solvers proposed in this paper. Table 2 displays the results of the numerical experiments for this test problem. The columns of Table 2 are essentially the same as the columns of Table 1. The first exception is the second column, where the values of $c$ are displayed rather than the order of approximation of the solver. An additional column was added, it shows the CPU gain factor (T-ratio) compared to the iterative solver used in  \cite{Turkel} .

\begin{table}[H]
\centering
\captionsetup{justification=centering}
\begin{tabular}{|r | r | r | r | r | r | r |}
\hline
     CPU  & c & Type &$N$ & $ \#$ iter.  & Time & T-ratio  \\  \hline
SM clust.  & 10 & iter.&229 &200& 122  & N/A \\ \hline
X5690 & 10&dir.& 197 & 1& 1.87  & 65 \\ \hline
iMac i7 &10&dir.& 197 &1& 1.63 & 75  \\ \hline
SM clust.  & 50& iter.&266 &280&289 & N/A \\ \hline
X5690 & 50&dir.& 280 & 1& 6.8 & 43  \\ \hline
iMac i7 &50&dir.& 280 &1& 6.9 & 42 \\ \hline
SM clust.  & 70& iter.&312 &893& 642 & N/A \\ \hline
X5690 & 70&dir.& 326 & 1& 12 & 54 \\ \hline
iMac i7 &70&dir.& 326 &1& 11 & 58 \\ \hline
SM clust.  & 80& iter.& $>402$ & N/A & N/A & N/A \\ \hline
X5690 & 80&dir.& 356 & 1& 12  & N/A \\ \hline
iMac i7 &80&dir.& 356 &1& 21  & N/A \\ \hline
\end{tabular}
\caption{Comparison of Direct and Iterative Solvers on the Second 3D Test Problem}
\end{table}

The range of parameters of the last test problems is more closely related to the realistic scenarios of the subsurface scattering problems. We can observe from the table that the direct solvers provide 42-75 times faster alternatives than the used iterative approach in \cite{Turkel}. They also allow the use of significantly larger grid sizes in comparison with the mentioned iterative solver on similar hardware. Table 2 indicates that the iterative solver could not find solutions on the grids greater than $402^3$. However, the sixth-order direct solver proposed in this paper was successfully applied to the last problem with $c=80$ and produced a solution satisfying the desired goal on both iMac and Xeon server. It must be noted that on the grids with the size greater than $356^3$, the significant advantage of the Xeon server in RAM becomes crucial for a rapid solution of the problem.  Remarkably, all calculations with direct solvers were conducted on a single Intel Core i7, 2.93 GHz processor with 16 Gb of RAM or on an Intel Xeon X5690 processor running at 3.47 GHz with 144 Gb of RAM, the frequency of which is similar to only one node in the Supermicro cluster consisting of 12 such nodes on which the iterative solutions were achieved in \cite{Turkel}. In the majority of our experiments, the sixth-order solver allows the achievement of the desired accuracy in less CPU time than the fourth-order direct solver, but the possibility of different grid steps in $x-$, $y-$, and $z-$ directions makes the fourth-order compact scheme an attractive alternative in some situations. Next, we will consider the numerical experiments in which parallel implementation of the developed direct algorithms was investigated.

\subsection{Scalability of the Proposed Direct FFT Solvers}

The results of the sequential implementation of the developed high-resolution direct methods presented in the previous subsection demonstrated higher efficiency of the developed approach in comparison to one of the best general iterative methods applied to the series of test problems. 

In this section, the scalability properties of the developed algorithms and their limitations are discussed. One of the advantages of the proposed methodology is the natural parallelization of the presented methods. 

This property, in the case of the second-order approximation compact finite-difference schemes and the trilinear finite element discretization, was discussed in several publications (see e.g. \cite{eo}). However, to the best of our knowledge, the detailed investigation of the scalability of the proposed high-resolution approach has not been considered. In the following subsections, we consider the solution of the 3D Helmholtz equation on the grid sizes up to $4096^3$ and the solution of the 3D convection-diffusion equation to demonstrate the robustness of the presented approach. 

\subsubsection{Helmholtz Equation with Constant-Coefficient}

The first series of parallel experiments used constant-coefficient $k^2$ with $a=20$, $b=0$, $c=10$, $\gamma=16$ and $\beta=12$. In these tests, a $500^3$ rectangular grid for the second (\ref{scheme2}), fourth (\ref{scheme4}) and sixth-order (\ref{scheme6}) compact schemes was utilized. Table 3 displays the results of the OpenMP implementation of the developed algorithms on a standard quad-core desktop. The solution time required for the parallel implementation of every considered compact scheme demonstrates near linear scalability. In all three cases, the solution wall time (in seconds) decreases by a factor of close to 2 as the number of OpenMP threads doubles.
\begin{table}[H]
\centering
\begin{tabular}{| c | c | c | c |}
\hline
order\textbackslash \# of threads&1&2&4\\\hline
$2^\text{nd}$&36.81 sec&18.77 sec &9.62 sec\\\hline
$4^\text{th}$&36.02 sec&16.71 sec&9.96 sec\\\hline
$6^\text{th}$&36.45 sec&17.62 sec &9.51 sec\\\hline
\end{tabular}
\caption{Desktop Solution Time for OpenMP}
\end{table}

To compare the performance of OpenMP and MPI parallelization of the developed direct solver, both implementations were run on a single node on Cori. In these numerical experiments, only the sixth-order algorithm on a $500^3$ rectangular grid was tested. The results of this comparison are presented in Figure 3.

Similar to the first test, the total computational time in both implementations are reduced by approximately half as the number of processing units is doubled. The parallelization becomes less effective for larger numbers of OpenMP threads or MPI processes as the benefit of splitting the tasks across multiple processing units decreases while overhead becomes a dominant component in the solution wall time. The results also demonstrate that the OpenMP implementation has a slightly better performance on a single node than MPI parallelization. This can be explained by the required communication between MPI processes.

\begin{figure}[H]
\centering
\captionsetup{justification=centering}
\caption{Computation Time OpenMP vs MPI\\ Restricted to a Single Node, $k^2 = const$}
\includegraphics[trim=1cm 2cm -1cm 0,width=0.5\textwidth]{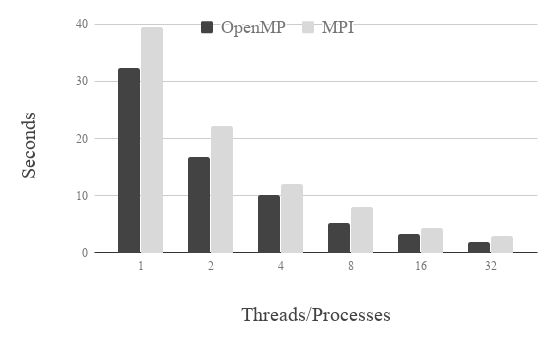}
\label{omp_vs_mpi_const}
\end{figure}

\subsubsection{Helmholtz Equation with Variable Coefficient}

Next, we consider the performance of the proposed parallel algorithms in the case of nonconstant-coefficient $k^2(z)$ with $a=10$, $b=9$, $c=10$, $\gamma = 9$ and $\beta = 10$. As in the previous example, the OpenMP and MPI implementations of the developed direct method were used to obtain approximate solutions of the problem (\ref{problem}) with the boundary conditions (\ref{bc}) on a sequence of grids. Both implementations give consistent results on all considered grids and different sets of the processing units. These results are almost identical to the convergence results of the sequential algorithm on corresponding grids. Tables \ref{conv_sec_table}, \ref{conv_fourth_table} and \ref{conv_sixth_table} demonstrate the convergence of the second, fourth and sixth-order sequential algorithms respectively. 

\begin{table}[H]
\centering
\begin{tabular}{|c|c|c|c|}\hline
&\textbf{max-err}&$L_2$\textbf{-err}&$L_2$\textbf{-res}\\ \hline 
$125^3$&5.7570466e-03&6.4986713e-03&4.7269292e-13 \\ \hline
$250^3$&
1.4853854e-03&
1.6510028e-03&
2.5846930e-12
\\ \hline
$500^3$&
3.7448165e-04&
4.1516358e-04&
6.5883688e-12
\\ \hline
\end{tabular}
\caption{Second-Order Convergence}
\label{conv_sec_table}
\end{table}

\begin{table}[H]
\centering
\begin{tabular}{|c|c|c|c|}\hline
&\textbf{max-err}&$L_2$\textbf{-err}&$L_2$\textbf{-res}\\ \hline 
$125^3$&
3.4493268e-05& 
3.5925614e-05&
3.6301725e-13 
\\ \hline
$250^3$&
2.1782070e-06&
2.2582699e-06&
1.9857221e-12
\\ \hline
$500^3$&
1.3726414e-07&
1.4187594e-07&
5.0832056e-12
\\ \hline
\end{tabular}
\caption{Fourth-Order Convergence}
\label{conv_fourth_table}
\end{table}

\begin{table}[H]
\centering
\begin{tabular}{|c|c|c|c|}\hline
&\textbf{max-err}&$L_2$\textbf{-err}&$L_2$\textbf{-res}\\ \hline 
$125^3$&
2.1875397e-06&
1.9909214e-06&
3.1581209e-13 
\\ \hline
$250^3$&
3.4942928e-08&
3.1643311e-08&
2.0541112e-12
\\ \hline
$500^3$&
5.5211108e-10&
4.9939925e-10&
5.2147803e-12 
\\ \hline
\end{tabular}
\caption{Sixth-Order Convergence}
\label{conv_sixth_table}
\end{table}

\begin{figure}[H]
\centering
\captionsetup{justification=centering}
\caption{Computation Time OpenMP vs MPI\\Restricted to a Single Node, $k^2 = k^2(z)$}
\includegraphics[trim=1cm 2cm -1cm 0,width=0.5\textwidth]{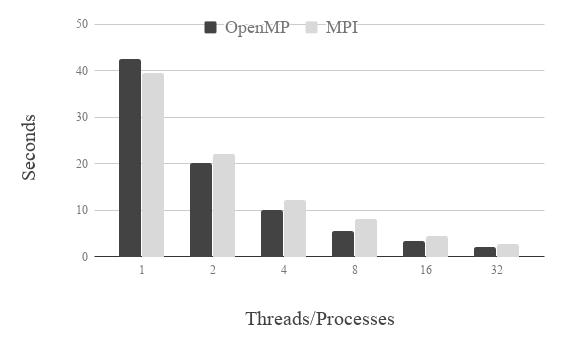}
\label{omp_vs_mpi_variable}
\end{figure}

The presented outcomes of the numerical experiments confirm the declared rate of convergence of the corresponding approximate solutions. To consider the comparison of the OpenMP and MPI implementations in the case of the variable problem coefficient, we use the same $500^3$ grid as in the constant-coefficient case. The sets of OpenMP threads and MPI processes were also chosen to be the same. Figure \ref{omp_vs_mpi_variable} shows nearly identical results as Figure \ref{omp_vs_mpi_const}. 

The limitations of OpenMP were observed in an attempt to run an experiment with the grid size of $1024^3$. The machines tested, including a single node on Cori, were unable to run this experiment due to a lack of memory. The experiment was repeated with the MPI implementation on one, two and four nodes on Cori. The attempts with both one and two nodes failed, again due to a lack of memory. However, four nodes successfully ran the program demonstrating the power of the MPI implementation. 

\begin{figure}[H]
\centering
\captionsetup{justification=centering}
\caption{Setup, Communication and Computation Time for MPI}
\includegraphics[trim=1cm 2cm -1cm 0,width=0.5\textwidth]{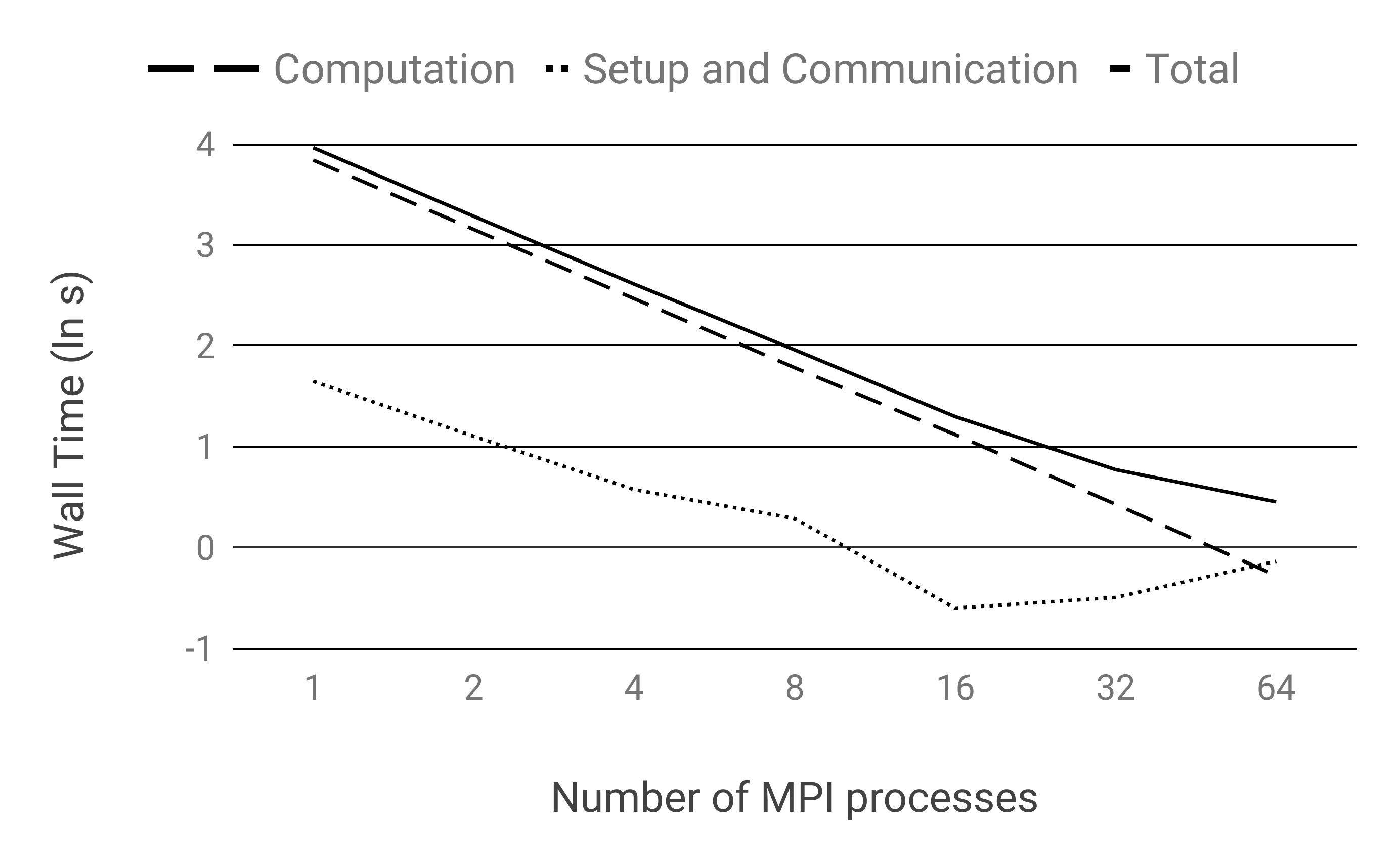}
\label{comm_lbnl}
\end{figure}

An experiment was conducted to further investigate the performance of the MPI parallelization. The setup, communication, and computation times were recorded as the number of MPI processes increased. Computation time refers to the time taken for all the MPI processes to complete the forward and inverse 2D DST, and tridiagonal solver steps. The communication time measures the longest time taken for the MPI processes to scatter the data to the appropriate processes and assigning data to a local array, while the setup time gives the time required to prepare the parallel environment. This test is run on a grid size of $512^3$ and Figure~\ref{comm_lbnl} shows the results.

We use the natural logarithm scale $\ln s$ for the vertical time-axis. On this graph, the computation time decreases almost linearly with the slope close to $-0.5$. On the other hand, the setup and communication time is decreasing only on the interval from 1 through 16 processes. It represents a small fraction of the total solution time on this interval. Beyond $16$ MPI processes, the setup and communication time is seen to increase. At $64$ MPI processes, the setup and communication time has exceeded the computation time in this particular test. This test demonstrated the limitation of the MPI implementation due to setup and communication times.  

To reduce the communication time the hybrid approach was considered. In this approach, one MPI process is assigned to one node. Further parallelization is achieved by using OpenMP threads available on every node. To present the performance of the hybrid implementation, a sequence of $2^k, k = 0, \ldots, 5$ MPI processes, i.e.\,nodes, and the matching sequence of OpenMP threads on every node were considered. Table \ref{hybrid_table} shows the computation times in seconds for each run on the same grid size $512^3$. Here, the number of OpenMP threads changes horizontally, and the MPI processes change vertically.
 
\begin{table}[H]
\centering
\begin{tabular}{|c|c|c|c|c|c|c|}\hline
n\textbackslash t&1&2&4&8&16&32\\ \hline
1&
39.61& 
21.73&
12.93&
8.83&
6.84&
5.62
\\ \hline
2&
19.87&
10.94&
6.48&
4.56&
3.48&
2.92
\\ \hline
4&
9.99&
5.66&
3.48&
2.52&
2.10&
1.82
\\ \hline
8&
5.27&
3.11&
1.99&
1.55&
1.37&
1.27
\\ \hline
16&
2.49&
1.42&
0.85&
0.64&
0.58&
0.58
\\ \hline
32&
1.85&
1.33&
1.09&
0.93&
0.76&
0.80
\\ \hline
\end{tabular}
\caption{Hybrid Implementation}
\label{hybrid_table}
\end{table}

The performance with strictly MPI, i.e.\,one thread, is similar to that shown in Figure \ref{omp_vs_mpi_variable}. However, in the case of restricting the hybrid algorithm to one node while increasing the number of OpenMP threads, the direct solver exhibits slightly different behavior in comparison with the results presented in Figure 4. This discrepancy can be explained by different utilization of the available OpenMP threads on every node. 

In the OpenMP implementation (Algorithm 1), the parallelization of the DST and the inverse DST steps is accomplished by the direct partition of the computational domain into a set of subdomains by horizontal $xy-$ planes. This method is justified in the previously run experiments since the number of grid points in the $z-$ direction was never below 125, while the number of cores available was always less than 32 in all CPUs under consideration. However, in the case of the hybrid approach, the number of grid points in $z-$ direction in the local MPI-subdomain is $N_z/N_{MPI}$, where $N_{MPI}$ is the number of MPI processes. This number becomes less than the number of OpenMP threads available in some of our numerical tests.

To balance the computational load between all available OpenMP threads in these situations, a different method of parallelization was adopted. As presented in the description of Algorithm 3, the sets of 1D FFTs on each of two stages of the 2D DST and inverse DST were divided into a sequence of subsets each of which was implemented using the multithread 1D FFTW application. The number of 1D FFTs in every subset, except the last one, is always the same as the number of OpenMP threads available. This method alleviates the parallelization restriction from $N_z/N_{MPI} > N_{OpenMP}$ to $N_z > N_{nodes}$, $N_x > N_{OpenMP}$, and $N_y > N_{OpenMP}$, where $N_{OpenMP}$ is the number of OpenMP threads available on a single node, and $N_{nodes}$ is the number of nodes used in the experiment. These restrictions are related only to the transform steps of the algorithm. Overall, this hybrid implementation produces similar results to the pure MPI parallelization on medium size grids, but it exhibits significantly better performance on the relatively large grids. This will be illustrated in our next experiments.

Further experiments were run to compare the performance of the MPI and hybrid implementations on a sequence of larger grids ranging from $512^3$ to $4096^3$. This was done to test the theory that reducing the number of MPI processes while maintaining the number of physical processors utilized will improve the computation time over MPI. Tables \ref{hybrid_vs_mpi_table_mpi} and \ref{hybrid_vs_mpi_table_hybrid} give these results.

\begin{table}[H]
\centering
\begin{tabular}{|c|c|c|c|c|c|c|}\hline
Grid & Nodes & MPI processes & Seconds\\ \hline
$512^3$&1&32&2.830525\\ \hline
$1024^3$&4&128&8.759851\\ \hline
$2048^3$&32&1024&40.465395\\ \hline
$4096^3$&256&4096&445.803343\\ \hline
\end{tabular}
\caption{Large Grid MPI}
\label{hybrid_vs_mpi_table_mpi}
\end{table}

\begin{table}[H]
\centering
\begin{tabular}{|c|c|c|c|c|c|c|}\hline
Grid & Nodes  & Processors& Seconds\\ \hline
$512^3$&1&32&7.793963\\ \hline
$1024^3$&4&128&16.911352\\ \hline
$2048^3$&32&1024&19.417831\\ \hline
$4096^3$&256&8192&27.522366\\ \hline
\end{tabular}
\caption{Large Grid Hybrid}
\label{hybrid_vs_mpi_table_hybrid}
\end{table}

While running this algorithm on a smaller grid size with a relatively small number of MPI processes, the MPI algorithm performed better than the hybrid implementation. This is likely due to the overhead required by OpenMP. When working with a grid size of $2048^3$ the communication time becomes a bottleneck in the MPI implementation since a total of $32$ nodes are required, which gives a total of $1024$ MPI processes. In the hybrid parallelization, only $32$ MPI processes are needed for the same grid size. This significantly reduces the communication time. The grid size $4096^3$ also demonstrates another limitation of the MPI algorithm. As mentioned in Section~\ref{sec: MPI} the MPI implementation is only able to utilize at most $4096$ processors. The hybrid implementation, however, significantly alleviate this restriction. Therefore, the hybrid implementation outperforms the strictly MPI program by $16$ times.

\subsubsection{Convection-Diffusion Equation}

This subsection presents the application of the developed direct parallel algorithms to the convection-diffusion equation  (\ref{cd_eq1}). Since this is simply an illustration of the diversity of applications of the proposed method, we restrict the consideration to the OpenMP implementation. The test problem under consideration can be presented as 
\begin{equation}
\nabla ^{2}u+\gamma \frac{\partial u}{\partial z} =0 , \ \  \textnormal{in}  \ \Omega, \label{cd_problem}
\end{equation}
where $ \gamma = -100 $ and on $\Omega = \left[0,\sqrt{2}\right]\times \left[0,\sqrt{2}\right] \times \left[0,1\right]$. With boundary conditions:
$u(x,y,0) = \sin\left(\frac{\pi x}{\sqrt{2}}\right)\sin\left(\frac{\pi y}{\sqrt{2}}\right), \ \ \ \ \ \ \ \ \ \ \ \ \ $
$u(0,y,z) = u(\sqrt{2},y,z) = u(x,0,z) = u(x,\sqrt{2},z) = 0, \ \ \ \ \ \ \ \ \ $
$u(x,y,1) = 2\sin\left(\frac{\pi x}{\sqrt{2}}\right)\sin\left(\frac{\pi y}{\sqrt{2}}\right), \ $ where $  \ 0 \le x,y \le \sqrt{2} $ and $\ \ \ \ \ $$0 \le z \le 1$. 
The analytic solution of the problem is given by
 \begin{align*}u = \sin\left(\frac{\pi x}{\sqrt{2}}\right)\sin\left(\frac{\pi y}{\sqrt{2}}\right)e^{-\gamma z/2}\frac{2e^{\gamma/2}\sinh(\sigma z)+\sinh(\sigma(1-z))}{\sinh\sigma}
 \end{align*}
 where $\sigma = \sqrt{\pi^2+ \gamma^2/4}$. The fourth-order convergence of the approximate solution to the analytic solution on a sequence of grids is presented in the following table. 
\begin{table}[H]
\centering
\begin{tabular}{|c|c|c|c|}\hline
 &	max-err &$L_2$-err & $L_2$-res \\ \hline
$64^3$ & 3.2612907e-03  & 4.6813690e-04&  1.5435312e-15 \\ \hline
$128^3$ & 2.0579387e-04  & 2.9792890e-05 &  4.6094565e-15  \\ \hline
$256^3$ & 1.2939970e-05  & 1.8507601e-06  &  9.9002420e-15  \\ \hline
$512^3$ & 8.1975702e-07  & 1.1559163e-07  &  3.2599029e-14  \\ \hline
\end{tabular}
\caption{Fourth-Order Convergence}
\label{cd_con}
\end{table}
Table \ref{CD_OMP_TABLE} gives the results of the parallel calculations on a single Cori node. The number of OpenMP threads in this test varies from one to eight. The table presents the wall time of the direct solver for two grid sizes, $256^3$ and $512^3$, using fourth-order approximation. We can see that the proposed parallel algorithm gives approximately a four times speed up from one thread, i.e. running sequentially, to eight threads. This test confirms the high efficiency and versatility of the developed parallel solver.
\begin{table}[H]
\centering
\begin{tabular}{| c | c | c |}\hline
OpenMP Threads &	$256^3$ &$512^3$ \\ \hline
1  &	7.858606 &	40.0879 \\ \hline
2 &	4.557636 &	24.553443 \\ \hline
4 &	2.666061 &	15.109893 \\ \hline
8 &	1.828526 &	10.878381 \\ \hline
\end{tabular}
\caption{Seconds to Compute}
\label{CD_OMP_TABLE}
\end{table}

\section{Conclusion}
In this paper, a direct parallel generalized FFT type algorithm was developed for a class of compact numerical approximations on a rectangular grid. The target applications of high-order compact approximation of 3D Helmholtz and convection-diffusion equations were considered on a sequence of grids. The developed algorithms represent highly accurate and scalable methods for the solution of the considered problems. The results demonstrated the efficiency of the OpenMP, MPI and hybrid implementations of the proposed parallel algorithms. This includes the vast improvement in computation time and the ability to apply these methods to other schemes with similar 3D stencils.
\section{Acknowledgements}
The authors gratefully acknowledge the financial support from the Sustainable Horizons Institute and Lawrence Berkeley National Laboratory(LBNL) in the form of the summer fellowships for all three authors at LBNL during summer 2018. The authors also thank the anonymous reviewers for their useful comments.

\section*{References}

\bibliographystyle{elsarticle-num}
\bibliography{bibliography}

\end{document}